\documentclass[11pt,twoside]{article}

\usepackage{epsfig,amsfonts,color}
\usepackage{amsmath, bbm, mathabx}
\usepackage[normalem]{ulem}

\usepackage{amssymb, palatino, geometry,url}
\usepackage{algorithmic}
\usepackage[noresetcount,lined,boxed]{algorithm2e} 

\usepackage[colorlinks=true,linkcolor=blue,citecolor=blue,urlcolor=blue]{hyperref}
\usepackage{subcaption}
\usepackage{amsthm}
\usepackage[utf8x]{inputenc}
\usepackage{textgreek}

\usepackage{mathrsfs}
\usepackage{multirow}

\usepackage{fancyhdr}
\usepackage{lineno}
\usepackage{float}
\usepackage[section]{placeins}
\usepackage[table]{xcolor}

\title{Random Field Optimization}
\author{Joshua L. Pulsipher, Benjamin R. Davidson, and Victor M. Zavala\thanks{Corresponding Author: victor.zavala@wisc.edu}\\
	{\small Department of Chemical and Biological Engineering}\\
	{\small \;University of Wisconsin-Madison, 1415 Engineering Dr, Madison, WI 53706, USA}}
\date{}

\begin{document}
	
\maketitle

\begin{abstract}
We present a new modeling paradigm for optimization that we call random field optimization. Random fields are a powerful modeling abstraction that aims to capture the behavior of random variables that live on infinite-dimensional spaces (e.g., space and time) such as stochastic processes (e.g., time series, Gaussian processes, and Markov processes), random matrices, and random spatial fields. This paradigm involves sophisticated mathematical objects (e.g., stochastic differential equations and space-time kernel functions) and has been widely used in neuroscience, geoscience, physics, civil engineering, and computer graphics. Despite of this, however, random fields have seen limited use in optimization; specifically, existing optimization paradigms that involve uncertainty (e.g., stochastic programming and robust optimization) mostly focus on the use of finite random variables. This trend is rapidly changing with the advent of statistical optimization (e.g., Bayesian optimization) and multi-scale optimization (e.g., integration of  molecular sciences and process engineering). Our work extends a recently-proposed abstraction for infinite-dimensional optimization problems by capturing more general uncertainty representations. Moreover, we discuss solution paradigms for this new class of problems based on finite transformations and sampling, and identify open questions and challenges. 
\end{abstract}

\noindent{\bf Keywords:} random fields; space-time; uncertainty; stochastic programming; infinite-dimensional optimization

\section{Introduction} \label{sec:intro}

Many complex systems are modeled using variables that live on continuous domains (e.g., space and/or time), making such models infinite-dimensional (i.e., variables are functions). These systems are difficult to model and optimize since they involve sophisticated mathematical objects such as partial differential equations (PDEs), derivatives, and integrals. Dynamic optimization and PDE-constrained optimization are common optimization modeling paradigms that involve such objects and these paradigms can be more broadly classified as infinite-dimensional optimization (InfiniteOpt) problems \cite{pulsipher2021unifying}. These problem classes arise in application domains such as dynamic system identification \cite{shin2019scalable}, model predictive control \cite{rawlings2000tutorial}, process intensification \cite{baldea2018dynamic}, and path planning \cite{ahmadzadeh2009multi}.  Specific applications include autonomous vehicles \cite{dadkhah2012survey}, process systems \cite{dimitriadis1995flexibility}, diffusive processes \cite{lu2021image}, microbial communities \cite{shin2019scalable}, structural design \cite{chen2010level}, molecular dynamics \cite{rapaport2004art}, and reaction systems \cite{graham2013modeling}. 
\\

The complexity of infinite-dimensional models is significantly exacerbated when these also capture random phenomena. Such phenomena can be modeled using random variables that live in finite-dimensional spaces (realizations live on finite/discrete domains such as vectors and matrices) or that live infinite-dimensional spaces (realizations are functions that live on continuous domains). Examples of finite-dimensional random variables are kinetic parameters or failure times, while examples of infinite-dimensional random variables include environmental disturbances (e.g., wind and temperature), uncertain transport fields (e.g., porosity), forecasts (e.g., time series), and molecular-level effects (e.g., random trajectories of particles) \cite{neckel2013random}. Incorporating infinite-dimensional random variables is key in engineering and scientific applications, but this practice has been quite limited in the optimization domain. One can attribute this to the fact that infinite-dimensional random variables involve complex mathematical objects such as stochastic differential equations and space-time kernel functions that are difficult to represent and handle computationally. For instance, stochastic processes \cite{ross1996stochastic} have been widely used to model random phenomena physics \cite{van1992stochastic}, biology \cite{allen2010introduction}, chemistry \cite{van1992stochastic}, neuroscience \cite{laing2010stochastic}, image processing \cite{geman1984stochastic}, and engineering \cite{graham2013modeling}. Classical models of stochastic processes are Wiener processes (i.e., Brownian motion), Markov processes, and Poisson processes. Stochastic processes and their impact on physical systems can be modeled using stochastic differential equations (SDEs) \cite{soong1973random}. SDEs capture the evolution of a stochastic variable, whose behavior is characterized by the time evolution of a probability density function. SDEs are a special case where a classical differential equation is perturbed by a random noise term (often a Wiener process) that models the randomness introduced in a dynamical system \cite{protter2005stochastic}. For instance, the diffusion of particles across a surface can be modeled via an SDE using a Wiener noise process and its solution (in terms of the resulting probability density function) yields the so-called Fokker-Planck equations for transient diffusion on a planar surface \cite{graham2013modeling}. Similarly, stochastic chemical kinetics (i.e., reactions occurring at low concentrations such that random molecular behavior is prominent) can be modeled as SDEs in which the extents of reaction are  characterized as Poisson processes; the solution of this system of SDEs for a general reaction network yields the so-called chemical master equations (i.e., forward Kolmogorov equation) \cite{graham2013modeling}. These examples help illustrate the significance of using stochastic models to describe dynamical systems. A key caveat of SDE models is that the underlying stochastic process does not have a closed-form solution and thus requires the use of stochastic calculus and/or simulation techniques to obtain the system evolution \cite{protter2005stochastic}. There is a variety of software tools available to accomplish this task such as \texttt{DifferentialEquations.jl} \cite{rackauckas2017differentialequations}, \texttt{SDE Toolbox} \cite{picchini2007sde}, and \texttt{SDE-MATH} \cite{kuznetsov2021sde}. We also note that the scope of stochastic processes and SDEs is generally limited to systems that evolve over time domains (as opposed to space-time domains). 
\\

Multi-stage stochastic optimization is a widely-used modeling paradigm to make decisions under uncertainties that evolve over time. Such formulations employ discrete-time stochastic processes (finite random variables) to model uncertainty; in this formulation, each stage seeks to optimize a random cost function that is conditional to the observed information from previous stages \cite{shapiro2021lectures}. The progressive gathering of information is captured in the form of a  scenario tree;  unfortunately, scenario trees grow exponentially with the number of stages \cite{shapiro2005complexity}. Diverse schemes have been developed to circumvent this complexity, such as the so-called stochastic dual dynamic programming (SDDP) algorithm  \cite{shapiro2011analysis}, which has been implemented in software tools such as \texttt{SDDP.jl} \cite{dowson2021sddp}. Multi-stage stochastic programming provides a natural framework to capture time-dependent uncertainties but the representation of such uncertainties is inherently finite-dimensional (over discrete times). Key questions that arise in this context are: how can we model stochastic programming problems that involve uncertainties that evolve over continuous time? Moreover, how can we characterize uncertainty that evolves over other domains (e.g., space and space-time)?
\\

Random field theory provides a powerful mathematical framework for characterizing uncertain factors that evolve over continuous (infinite) domains such as space-time. Under this paradigm, uncertainty is modeled as a random function/manifold/field \cite{adler2007random}. Each realization of a random field is a deterministic function defined over the infinite domain and there is an underlying probability density that characterizes the probabilities of such realizations. This general representation captures a wide range of uncertainties such as stochastic processes (domain is continuous time), Gaussian processes (domain is space, time, or space-time), point processes (domain collapses to a finite set of points), and random vectors/matrices (a field defined over a finite set of points). Random field models have been widely studied in physics \cite{adler2007random}, neuroscience \cite{brett2003introduction}, computer graphics \cite{wang2013markov}, geoscience \cite{christakos2012random,harbaugh1968fourier}, and civil engineering \cite{chen2010level}. This rich theoretical base has enabled the development of advanced tools for modeling and analysis; for instance, random field models have been used for identifying abnormal behavior in brain imaging data (modeled as space-time random fields) \cite{adler2008some}. Moreover, random field theory provides tools that facilitate the characterization of extreme  behavior (e.g., excursion probabilities) and has fundamental connections with other fields of mathematics such as topology. Surprisingly, random fields have seen limited use in the context of optimization. To the best our best knowledge, random fields have only been used to model uncertainty of material properties in the context of topological optimization in civil engineering \cite{chen2010level}. Moreover, random fields (Gaussian processes) have recently been used in the context of Bayesian optimization for machine learning and to build surrogate models in model predictive control, but such work does not account for uncertainties that evolve over continuous space-time domains. 
\\

Recently, in \cite{pulsipher2021unifying}, we presented a unifying abstraction for InfiniteOpt problems. This abstraction utilized a general treatment of infinite domains (e.g., space, time, and/or uncertainty) and measure theory to model associated variables, objectives, and constraints. This enables us to commonly express diverse optimization formulations stemming from areas such as dynamic, PDE-constrained, stochastic, and semi-infinite optimization; which enables theoretical crossover and development of new formulation paradigms. In this context, we observed that this framework can facilitate the incorporation of random fields to model uncertainty that evolves over infinite domains. 

\begin{figure}[!htb]
	\includegraphics[width=\textwidth]{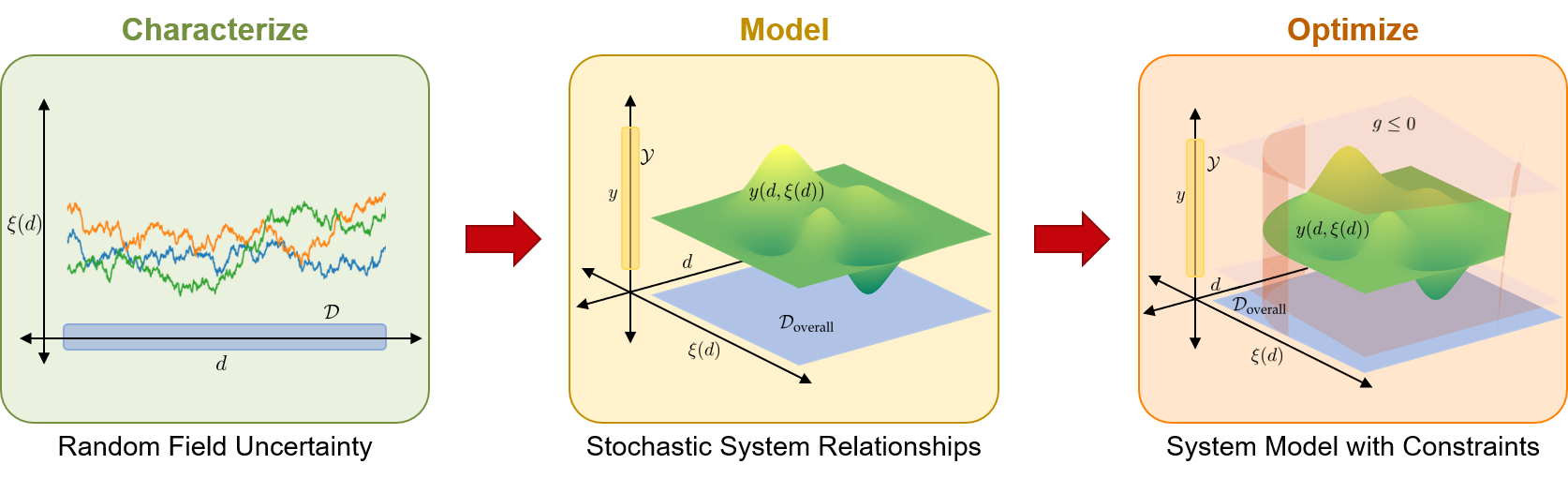}
	\centering
	\caption{Summary of proposed random field optimization framework.}
	\label{fig:abstract}
\end{figure}

In this work, we formalize the incorporation of random field theory into InfiniteOpt problems, establishing a new modeling paradigm that we call \emph{random field optimization} (RFO). This modeling paradigm aims to unify and facilitate the development of models that diverse sources of uncertainty. The framework characterizes such uncertainties via random fields and aims to integrate associated mathematical objects (e.g., SDEs and kernel functions) directly into optimization formulations. Figure \ref{fig:abstract} summarizes this approach. RFO provides an intuitive framework that will serve as the foundation for future research in modeling, analyzing, and solving this new general class of InfiniteOpt problems. Moreover, this abstraction can enable further development of software tools such as \texttt{InfiniteOpt.jl} to make the optimization of infinite-dimensional systems more accessible.
\\

The paper is structured as follows. Section \ref{sec:random_fields} introduces random fields and establishes mechanisms to characterize uncertainty. Section \ref{sec:modeling} discusses how random fields are incorporated into system models. Sections \ref{sec:optimization} describes how to model RFO problems. Section \ref{sec:cases} presents illustrative case studies and Section \ref{sec:conclusion} closes with final remarks.

\section{Random Fields} \label{sec:random_fields}

We detail relevant definitions, properties, characterizations and metrics from random field theory and establish basic notation. The material covered here is not intended to be a comprehensive review; we refer the reader to \cite{adler2008some,christakos2017spatiotemporal, chung2020introduction} for a thorough review of random field theory.

\subsection{Definition}

\begin{figure}[!htb]
	\includegraphics[width=0.4\textwidth]{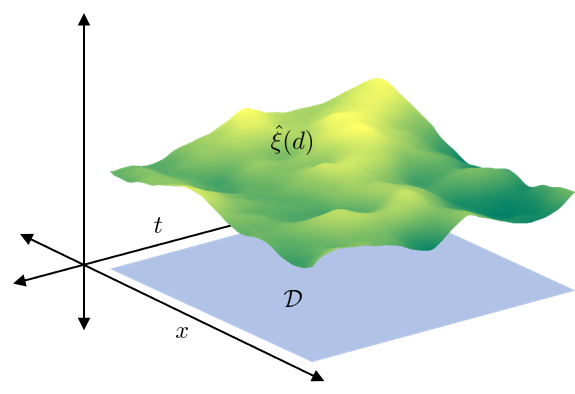}
	\centering
	\caption{Illustration of a realization $\hat{\xi}(d)$ of random field $\xi(d)$ where $n_\xi = 1$.}
	\label{fig:field_sample}
\end{figure}

Random field theory refers to the study of random functions defined over Euclidean spaces (e.g., space and/or time) \cite{adler2008some}. A random field is a random function $\xi : \mathcal{D} \mapsto \mathcal{D}_{\xi(d)}$ whose input domain $\mathcal{D} \in \mathbb{R}^{n}$ is a general topological space:
\begin{equation}
    \xi(d) \in \mathcal{D}_{\xi(d)}, \ d \in \mathcal{D}
    \label{eq:random_field}
\end{equation}
where $\mathcal{D}_{\xi(d)}$ is a function space that is the co-domain of the random field (i.e., its realizations are deterministic functions $\hat{\xi} : \mathcal{D} \mapsto \mathbb{R}^{n_\xi}$). Here, a realization of the random field $\hat{\xi}(d) \in \mathcal{D}_{\xi(d)}$ (also called a sample function) is a deterministic function, as illustrated in Figure \ref{fig:field_sample}.

Random field theory is a flexible paradigm that is amenable to diverse engineering applications such as petroleum reservoir management, diffusive surfaces, fluid dynamics, atmospheric modeling, granular and porous materials, stochastic chemical kinetics, and molecular dynamics. Modeling the height of the ocean surface over a certain spatial area provides an intuitive example of a random field. Here, a realization denotes a snapshot that describes a possible ocean surface with defined waves. Moreover, we can add a temporal dependency such that a realization describes how a set of waves propagate over the surface over time. A key observation here is that the surface height exhibits continuous correlation over space and time and such correlations can be captured using models proposed in random field theory. Continuity is key, as this enforces regularity and restricts the nature of the variability in the surface observed. 
\\

Following the terminology from the recently-proposed unifying abstraction for InfiniteOpt problems, we use $d$ to denote an infinite parameter (e.g., space and time) that indexes the infinite domain $\mathcal{D}$ (typically referred to as a topological manifold in the random field literature) \cite{pulsipher2021unifying,adler2007random}. Equivalently, we interpret a random field $\xi(d)$ as a set of random variables $\xi_d$ indexed over the domain $\mathcal{D}$:
\begin{equation}
    \{\xi_d : d \in \mathcal{D}\}.
    \label{eq:random_field_alt}
\end{equation}
Moreover, these are defined with a set of density functions $\{F_{d_1}, \dots, F_{d_n}\}$ that are in turn defined on the Borel space $\mathscr{B}^{nm}$ that satisfies:
\begin{equation}
	F_{d_1}, \dots, F_{d_n}[B] = P[(\xi(d_1),\dots \xi(d_n)) \in B], \ B \in \mathscr{B}^{nm}
	\label{eq:finite_dists}
\end{equation}
where $m$ is the dimension of $\mathcal{D}_{\xi(d)}$ and $n$ is an arbitrary positive integer. These are often referred to as the finite-dimensional densities of a random field and are central in establishing fundamental properties for a random field, such as those presented in Section \ref{sec:field_properties}.
\\

Following \eqref{eq:random_field_alt}, a random field can be interpreted as the {\em generalization of a multivariate random variable} where the indexing set is continuous (infinite-dimensional). For instance, a discrete time-series $\{\xi_t, t \in \mathcal{D}_t\}$ is a multivariate random variable and this can be generalized as a temporal random field when $\mathcal{D}_t$ becomes a continuous time interval. In this work, we will favor the use of \eqref{eq:random_field}, since function notation is more amenable to infinite-dimensional system models (which is the focus of this work). Note, however, that the term infinite-dimensional arises from the observation that system variables/parameters are indexed over a space $\mathcal{D}$ of infinite cardinality, as exemplified in \eqref{eq:random_field_alt}.

\subsection{Properties} \label{sec:field_properties}

We begin by establishing the concept of mean and covariance functions of a random field, which are basic constructs needed to derive additional random field properties. These can be thought of as continuous generalizations of mean vectors and covariance matrices used in random matrix theory. 
\\

The mean function of a random field $\mu(d)$ is a deterministic function that denotes the expected value over the random field manifold $\mathcal{D}$:
\begin{equation}
	\mu(d) := \mathbb{E}_{\xi(d)}[\xi(d)].
	\label{eq:mean_function}
\end{equation}
The covariance function $\Sigma : (\mathcal{D}, \mathcal{D}) \mapsto \mathbb{R}$ is a deterministic function and captures the correlation/variability exhibited over the random field manifold:
\begin{equation}
	\Sigma(d, d') := \mathbb{E}_{\xi(d)}\left[(\xi(d) - \mu(d))(\xi(d') - \mu(d'))\right].
	\label{eq:covar_function}
\end{equation}
We can see that, if $\mathbb{E}_{\xi(d)}[\xi(d)^2]$ is finite for all $d \in \mathcal{D}$, then $\Sigma(d, d')$ is finite for all $d,d' \in \mathcal{D}$; moreover, $\Sigma(d, d')$ is assumed to be a non-negative definite function \cite{adler2008some}. We will see in Section \ref{sec:field_characteriations} that the mean and covariance functions are often useful in defining (Gaussian) random fields. 
\\

A common consideration in modeling random fields is that of {\em homogeneity}. A random field is said to be homogeneous if its finite-dimensional densities are translationally invariant. Formally, this indicates that the joint density of any $k$ random variables $\{\xi(d_1), \dots, \xi(d_k)\}$ is equivalent to that of the random variables $\{\xi(d_1 + \tau), \dots, \xi(d_k + \tau)\}$ for arbitrary $\tau \in \mathbb{R}$. A homogeneous random field thus exhibits a constant mean function and its covariance function is only a function of the difference $d-d'$:
\begin{equation}
	\begin{aligned}
		\mu(d) &= \beta \\
		\Sigma(d, d') &= \Sigma(d - d')
	\end{aligned}
	\label{eq:homo_props}
\end{equation}
where $\beta \in \mathbb{R}^{n_\xi}$ is a constant. Subsequently, a random field is said to be weakly homogeneous or stationary if it only satisfies \eqref{eq:homo_props}.
\\

Another key consideration (that will be particularly relevant in Section \ref{sec:modeling}) pertains to the continuity and the differentiability of random fields. Continuity and differentiability are key in capturing behavior, because such properties naturally regularize correlation; for example, these properties indicate that the value of a random field at a space-time location is not arbitrarily far away from that at a neighboring point. There are a considerable number of definitions used in the literature to derive and discuss these properties. A common group of definitions, that are convenient for our purposes, define the continuity and differentiability of random fields with respect to their sample functions (realizations). 
\\

A random field is sample function continuous at a point $d^*$ if for every sequence $\{d_n\}$ that satisfies $||d_n - d^*|| \rightarrow 0$ as $n \rightarrow \infty$ we have that:
\begin{equation}
	P\left[||\xi(d_n) - \xi(d^*)|| \rightarrow 0 \;\; \text{as} \; n \rightarrow \infty \right] = 1.
	\label{eq:field_continuity}
\end{equation}
This essentially amounts to almost surely having sample functions that are continuous at a point $d^*$. Similarly, we say a random field is sample function differentiable at a point $d^*$ if we have that:
\begin{equation}
	P\left[\frac{\partial \xi(d)}{\partial d_\ell} \; \text{exists}\right] = 1
	\label{eq:field_differentiability}
\end{equation}
where $d_\ell$ is the $\ell^\text{th}$ element of $d$. Hence, we can define a random field as continuous and/or differentiable if these conditions hold almost surely for its (deterministic) sample functions. 

\subsection{Characterization} \label{sec:field_characteriations}

\begin{figure}[!htb]
	\includegraphics[width=0.4\textwidth]{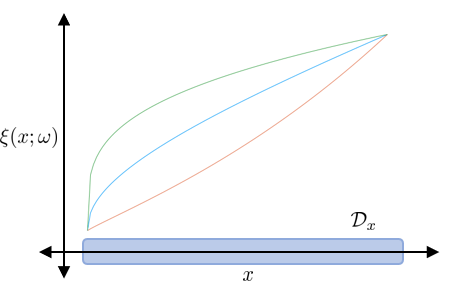}
	\centering
	\caption{Illustration of realizations $\hat{\xi}(x; \hat{\omega})$ from a simple random function $\xi(x; \omega)$ which is a random field over spatial position $x \in \mathcal{D}_x$.}
	\label{fig:analytic_field}
\end{figure}

Since random field theory provides a broad framework for random functions, there exist numerous ways to model/characterize them. For an intuitive initial step, let us consider a simple analytic function $\xi(d; \omega)$ with a parametric dependence on a random variable $\omega$. Figure \ref{fig:analytic_field} depicts such a function, where $\xi(x; \omega) = \exp(x^\omega)$ denotes a random field. In practice, however, such functions can be modeled directly with random variables and need not be treated rigorously as random fields. 

The most prevalent random field characterizations arising in applications are Gaussian random fields, which can be completely specified via their second-order statistics (i.e., mean and covariance functions). More precisely, a random field is Gaussian if all of its finite-dimensional density functions are multi-variate Gaussian (this generalizes Gaussian random variables over a general topological domain $\mathcal{D}$). This unique property enables a wealth of analysis and modeling techniques, since the mean and covariance are available in convenient analytical forms for multi-variate Gaussian distributions. 
\\

\begin{figure}[!htb]
     \centering
     \begin{subfigure}[b]{0.31\textwidth}
         \centering
         \includegraphics[width=\textwidth]{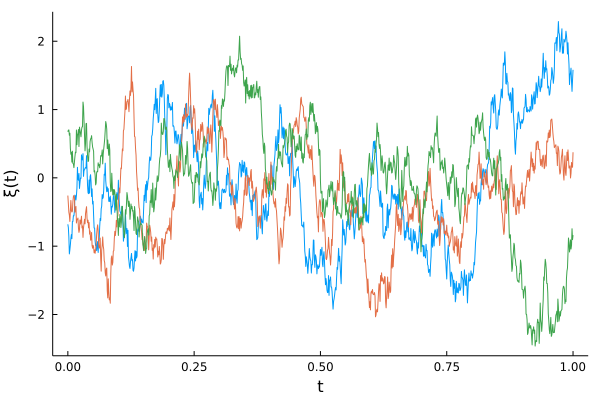}
         \caption{Ornstein-Uhlenbeck}
     \end{subfigure}
     \quad
     \begin{subfigure}[b]{0.31\textwidth}
         \centering
         \includegraphics[width=\textwidth]{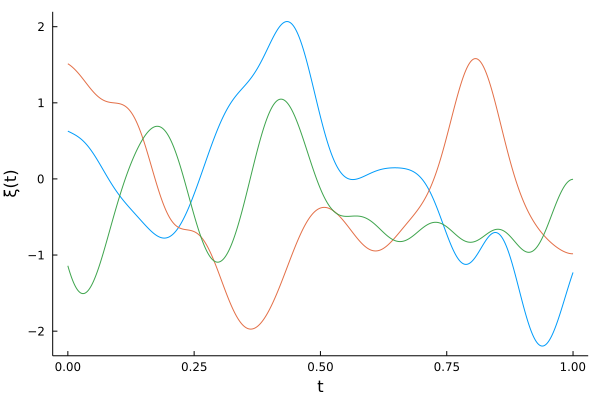}
         \caption{Squared Exponential}
     \end{subfigure}
     \quad
     \begin{subfigure}[b]{0.31\textwidth}
         \centering
         \includegraphics[width=\textwidth]{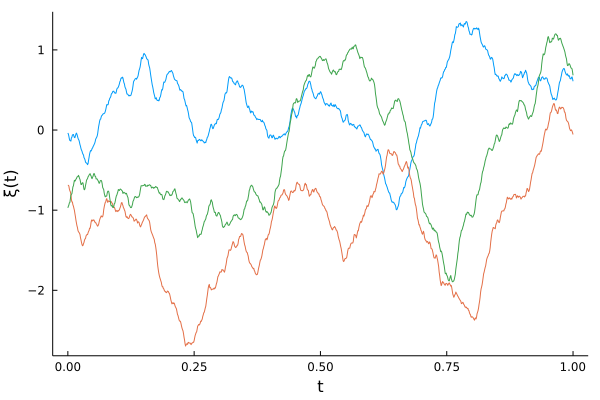}
         \caption{Matern}
     \end{subfigure}
    \caption{Realizations of Gaussian random fields with varied covariance functions; where appropriate we use $\beta = 0.1$ and $\nu,\sigma = 1$. Note how the smoothness of the realizations are affected by the form of the covariance function.}
    \label{fig:covar_kernels}
\end{figure}

Without a loss of generality, the mean function is often assumed to be zero, such that the field is fully characterized via the choice of the covariance function. The prevalence of Gaussian random fields/processes in a wide variety of disciple areas (e.g., data science, statistics, and engineering) has motivated extensive development of models of covariance functions (often referred to as kernel functions). Common covariance functions include linear $\Sigma_\text{L}(d, d')$, squared exponential $\Sigma_\text{SE}(d, d')$, Ornstein-Uhlenbeck $\Sigma_\text{OU}(d, d')$ (also called the exponential kernel), and Matern $\Sigma_\text{M}(d, d')$:
\begin{equation}
	\begin{aligned}
		\Sigma_\text{L}(d, d') &:= d^Td' \\
		\Sigma_\text{SE}(d, d') &:= \sigma \exp\left(- \frac{||d-d'||^2}{2\beta^2}\right) \\
		\Sigma_\text{OU}(d, d') &:= \sigma \exp\left(- \frac{||d-d'||}{\beta}\right) \\
		\Sigma_\text{M}(d, d') &:= \frac{2^{1-\nu}}{\Gamma(\nu)}\left(\frac{\sqrt{2\nu}||d-d'||}{\beta}\right)^\nu K_\nu\left(\frac{\sqrt{2\nu}||d-d'||}{\beta}\right)
	\end{aligned}
\end{equation}
where $\sigma,\beta \in \mathbb{R}$ are scalar parameters, $K_\nu$ is a modified Bessel function with order $\nu \in \mathbb{R}$, and $\Gamma(\nu)$ is the gamma function \cite{kang2017slope}. 
Selecting a particular covariance function amounts to enforcing a prior on the random field that will affect its fundamental behavior. For example, the Ornstein-Uhlenbeck function (which corresponds to Brownian motion) is nowhere differentiable, while the squared exponential function gives rise to smooth field realizations that are almost surely differentiable. Figure \ref{fig:covar_kernels} illustrates this for various covariance functions. In addition to differentiability, we observe how the choice of the covariance function influences other field properties; for instance, covariance functions that only depend on the difference $d-d'$ are stationary and those that only depend on the distance $||d-d'||$ are stationary and isotropic (i.e., they have no directional dependence). There exist diverse software implementations of Gaussian random fields such as \texttt{GaussianRandomFields.jl} in \texttt{Julia}, \texttt{RandomFields} in \texttt{R}, and \texttt{GenPak} in \texttt{C++} that make it straightforward to build Gaussian fields and obtain sample functions. We refer the reader to \cite{liu2019advances} for an in-depth review of modeling Gaussian random fields.

Non-Gaussian field models have also been proposed in the literature, which include $\chi^2-$, $t-$, and $F-$ random fields \cite{worsley1994local}. As in the finite-dimensional case, these are constructed via a linear combination of Gaussian random fields that are independently and identically distributed \cite{chung2020introduction}. For example, we can construct a $\chi^2-$ random field $\xi_{\chi^2}(d)$ via the linear combination:
\begin{equation}
	\xi_{\chi^2}(d) = \sum_{i \in \mathcal{I}}\xi_{\text{G}, i}(d)
	\label{eq:chi_squared_field}
\end{equation}
where $\xi_{\text{G}, i}(d), \ i \in \mathcal{I}$, are i.i.d. Gaussian random fields with a zero mean and unit variance. Such fields can help capture a wider range of random field behavior. Moreover, we can construct random fields (and their associated covariance functions) via statistical analysis of the ensemble solution to an underlying system of RDEs (or SDEs) to model the uncertainty of interest \cite{christakos2017spatiotemporal,rozanov2013random}. This typically requires the use of stochastic calculus, but there are tools such as \texttt{DifferentialEquations.jl} that can be used to obtain accurate numerical ensembles of RDE/SDEs \cite{rackauckas2017differentialequations}.

\subsection{Measures} \label{sec:fields_measures}

Random fields are complex objects that often need to be projected/reduced via measures in order to be analyzed and optimized for \cite{pulsipher2021unifying}. Specifically, a measure will summarize the random field into a deterministic function or quantity. For instance, the expectation operator $\mathbb{E}_{\xi(d)}$ shown in \eqref{eq:mean_function} is a measure that summarizes the random field into a deterministic mean function $f(d)$ that only depends on the domain $\mathcal{D}$. Moreover, we could employ more general measure operators $M_d$ (e.g., integral over the domain $\mathcal{D}$) to further reduce the mean function to a scalar value. Hence, the use of measure operators helps us develop summarizing statistics and quantities that help analyze the behavior of a random field. Such summarizing statistics are also key in formulating optimization problems, as such problems often require the representation of objectives and constraints as deterministic quantities. 
\\

One of the most prevalent measures used in random field theory is the so-called {\em excursion probability} (also known as the exceedance probability):
\begin{equation}
	P_{\xi(d)}\left\{\sup_{d \in \mathcal{D}} \xi(d) \geq u\right\}
	\label{eq:excursion_prob}
\end{equation}
where $u \in \mathbb{R}$ is a threshold parameter and we take $\xi(d)$ to be an $\mathbb{R}$-valued field (this can be extended to vector valued fields \cite{fournier2018geometric}). The measure \eqref{eq:excursion_prob} is at the focus of most of the random field literature \cite{adler2007random}. This is because this measure summarizes the excursion behavior of a random field and thus helps characterize extreme values (i.e., the peaks of a function), which is useful in many applications. For instance, the excursion probability is used for hypothesis testing in neuroscience, geoscience, civil engineering and other fields to determine significant deviations relative to a certain threshold. Other potential uses include quantifying the likelihood that a field surpasses a certain limit/tolerance, which has implications for quality control, computing escape probabilities and times (e.g., in molecular dynamics), and in extreme event quantification. 

\begin{figure}[!htb]
	\includegraphics[width=0.4\textwidth]{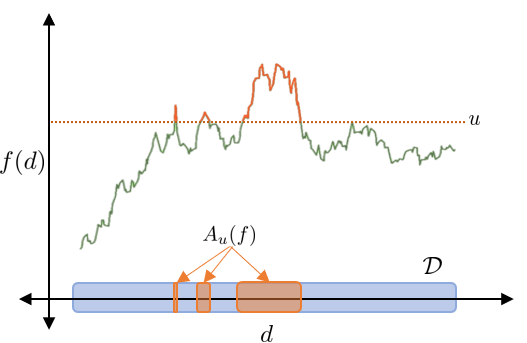}
	\centering
	\caption{Illustration of a 1D excursion set $A_u(f)$ on an arbitrary function $f(d)$.}
	\label{fig:excursion_set}
\end{figure}

To better understand \eqref{eq:excursion_prob} and discuss ways to compute it, it is necessary to define the notion of an excursion set $A_u$ for an arbitrary function $f$:
\begin{equation}
	A_u(f) := \{d \in \mathcal{D} : f(d) \geq u\}.
	\label{eq:excursion_set}
\end{equation}
The excursion set is the subdomain of $\mathcal{D}$ where $f$ exceeds $u$. Figure \ref{fig:excursion_set} illustrates such a set for a 1D function. Here, we observe that the excursion set of a 1D function can be characterized by a number of disjoint intervals (subdomains) in $\mathbb{R}$. This means that its topology can be described by the number of upcrossings that occur. This observation is the basis for extensive theoretical development for 1D random fields (e.g., stochastic processes over time) \cite{leadbetter2012extremes}. From this, a number of closed-form representations of \eqref{eq:excursion_prob} have been developed; for example, Rice's formula is commonly used  for computing the expected number of upcrossings (approximating \eqref{eq:excursion_prob}) for a zero-mean, stationary, continuous Gaussian field on $d \in [0, 1]$):
\begin{equation}
	\mathbb{E}[N_u] = \frac{\sqrt{\lambda}}{2\pi\sigma}\exp\left(\frac{-u^2}{2\sigma^2}\right)
	\label{eq:rice_formula}
\end{equation}
where $N_u$ denotes the number of upcrossings, $\sigma^2 = \mathbb{E}_{\xi(0)}[|\xi(d)|^2] = \Sigma(d)$ is the covariance function evaluated at $d = 0$, and $\lambda = -\Sigma''(0)$ is the second derivative of the covariance function evaluated at $d = 0$ \cite{adler2007random}. 
\\

The notion of upcrossings is ill-defined for (random) functions with $|\mathcal{D}| > 1$, making the analysis of high-dimensional excursion sets/probabilities significantly more challenging. However, in these cases, it is possible to use the expected Euler characteristic (EC) $\mathbb{E}[\mathcal{X}(A_u(\xi(d)))]$ as an approximation. It can be shown that, under certain regularity conditions, the error of this approximation decays exponentially with increased $u$; this makes the expected EC a useful construct used for  characterizing the behavior of random fields  \cite{adler2007random}. 
\\

The EC measures the topology/geometry of excursion sets (e.g., measuring the number of disjoint intervals in 1D deterministic sets and the number of connected components minus the number of holes for 2D deterministic sets) and capture the critical points (i.e., extrema) of a field since they define its topological features \cite{smith2021euler}. We refer the reader to \cite{smith2021euler} for an intuitive introduction to computing the EC of field excursions. 
\\

For certain types of random fields, there exist analytic representations of the expected EC characteristic \cite{adler2008some}; when these are not available, the expected EC can be computed via Monte Carlo (MC) sampling using the samples $\{\hat{\xi}_k(d) : k \in \mathcal{K}\}$:
\begin{equation}
	\mathbb{E}_{\xi(d)}[\mathcal{X}(A_u(\xi(d)))] \approx |\mathcal{K}|^{-1} \sum_{k \in \mathcal{K}} \mathcal{X}(A_u(\hat{\xi}_k(d))).
	\label{eq:mc_expected_ec}
\end{equation}
We can also approximate the excursion probability in \eqref{eq:excursion_prob} directly via MC sampling \cite{adler2012efficient} as:
\begin{equation}
	P_{\xi(d)}\left\{\sup_{d \in \mathcal{D}} \xi(d) \geq u\right\} = \mathbb{E}_{\xi(d)}\left[\mathbbm{1}\left(\sup_{d \in \mathcal{D}} \xi(d) \geq u\right)\right] \approx |\mathcal{K}|^{-1} \sum_{k \in \mathcal{K}} \mathbbm{1}\left(\sup_{d \in \mathcal{D}} \hat{\xi}_k(d) \geq u\right)
	\label{eq:mc_excursion_prob}
\end{equation}
where $\mathbbm{1}(\cdot)$ is the indicator function. There are a number of alternative measures (metrics) used in the random field literature, we refer the reader to \cite{adler2007random} for more details.

\section{Modeling with Random Fields} \label{sec:modeling}

In this section, we discuss how to incorporate random fields into systems models commonly encountered in applications. We also overview how the resulting stochastic models can be transformed into finite representations that are amenable for simulation and/or optimization. Much of this discussion includes a general overview of random and stochastic differential equation theory; we refer the interested reader to \cite{graham2013modeling, neckel2013random, soong1973random, cinlar2013introduction} for a rigorous treatment of these topics.

\subsection{Characterizations} \label{sec:model_characterize}

Once we have a random field representation of the uncertainty that affects a given system, we can incorporate it into the governing equations (e.g., differential and algebraic equations) to yield a stochastic model that we can simulate and/or use in an optimization formulation (which we discuss in Section \ref{sec:optimization}). We refer to governing equations that embed random field components as random differential-algebraic equations (RDAEs). A key observation that we start with is that a function that depends on a random field is itself a random field (i.e., $f(\xi(d))$ is a random field for arbitrary $f$); as such, as we construct algebraic and differential equations that incorporate a random field $\xi(d)$ they become inherently random in nature. The resulting stochastic model can be interpreted as a characterization for a high-dimensional random field model, which further highlights the connection between random field theory and random differential equations.

\subsubsection{Random Differential-Algebraic Equations}

\begin{figure}[!htb]
	\includegraphics[width=0.4\textwidth]{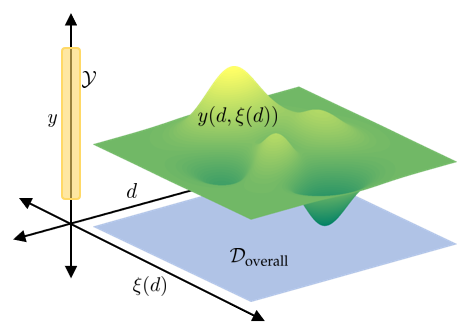}
	\centering
	\caption{High-level illustration of a random field variable $y(d, \xi(d))$.}
	\label{fig:variable}
\end{figure}

In the most general sense, we can define our model by directly embedding random field uncertainty into the DAEs that describe a system giving rise to a set of RDAEs. To do so, we consider system (state) variables $y : (\mathcal{D}, \mathcal{D}_{\xi(d)}) \mapsto \mathcal{Y} \subseteq \mathbb{R}^{n_y}$ which are functions over the system domain $\mathcal{D}$:
\begin{equation}
	y(d, \xi(d)) \in \mathcal{Y}
	\label{eq:model_vars}
\end{equation}
where we explicitly show their dependence on $\xi(d)$ to emphasize their inherent randomness. Figure \ref{fig:variable} depicts such a modeling variable. 

Following the notation proposed in \cite{pulsipher2021unifying}, we define the differential operator $D : \mathcal{D} \mapsto \mathcal{D}$ that operates on variables $y(d, \xi(d))$:
\begin{equation}
	Dy(d, \xi(d)) \in \mathbb{R}^{n_D}.
	\label{eq:model_derivs}
\end{equation}
This general operator captures partial derivatives $\partial y(d, \xi(d))/\partial d_\ell$, Laplace operators $\nabla^2y(d, \xi(d))$, and stochastic differentials $dy(d, \xi(d))$. With these objects we can now define general random DAEs of the form:
\begin{equation}
	g(Dy, y(d, \xi(d)), \xi(d), d) = 0
	\label{eq:random_daes}
\end{equation}
where $g(\cdot) \in \mathbb{R}^{n_g}$ are vector-valued functions. This form is general and captures a large collection of application systems. An important note is that classical differential operators like $\partial /\partial d_\ell$ can only be used if $\xi(d)$ is sample function differentiable in accordance with \eqref{eq:field_differentiability}. Otherwise, \eqref{eq:random_daes} needs to be carefully setup with SDEs derived from appropriate stochastic calculus constructs (e.g., It\^{o} or Stratonovich calculus) \cite{graham2013modeling,soong1973random}.
\\

We also note that concerted efforts have been recently made to develop theory for RDAEs (especially in the context of SDEs). However, these have been limited to a few prescribed forms thus far. Relevant recent works for further reading on this topic include \cite{winkler2006stochastic, suthar2021explicit}.

\subsubsection{Random Differential Equations}

Due to the complexity associated with stochastic differential equations, theory is not well-developed for general random systems and the literature commonly considers a set of RDEs of the form:
\begin{equation}
	\begin{gathered}
		Dy(d, \xi(d)) = f(y(d, \xi(d)), \xi(d), d) \\
		q(y(d_0, \xi(d_0)), \xi(d_0), d_0) = 0, \ d_0 \in \mathcal{D}_0
	\end{gathered}
	\label{eq:model_rdes}
\end{equation}
where $f: (\mathcal{Y}, \mathcal{D}_{\xi(d)}, \mathcal{D}) \mapsto \mathbb{R}^{n_d}$ is an arbitrary right-hand side function, $q : (\mathcal{Y}, \mathcal{D}_0) \mapsto \mathbb{R}$ encodes the initial/boundary conditions, and $\mathcal{D}_0$ is the set of initial and/or boundary points \cite{soong1973random, neckel2013random}. Three types of RDEs are encountered in the literature; these are distinguished by the way randomness (characterized by $\xi(d)$) enters the system. These are (in order of increasing difficulty to analyze and solve) random initial/boundary conditions, random inputs (via an inhomogeneous random term), and random coefficients \cite{soong1973random}. 
\\

\begin{figure}[!htb]
    \centering
	\begin{subfigure}[b]{0.31\textwidth}
        \centering
        \includegraphics[width=\textwidth]{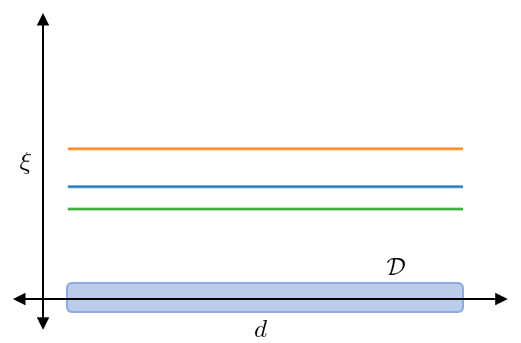}
        \caption{Static}
        \label{fig:static_field}
    \end{subfigure}
	\quad
    \begin{subfigure}[b]{0.31\textwidth}
        \centering
        \includegraphics[width=\textwidth]{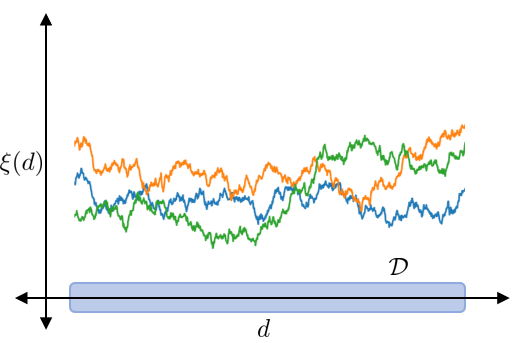}
        \caption{Correlated}
        \label{fig:correlated_field}
    \end{subfigure}
    \quad
    \begin{subfigure}[b]{0.31\textwidth}
        \centering
        \includegraphics[width=\textwidth]{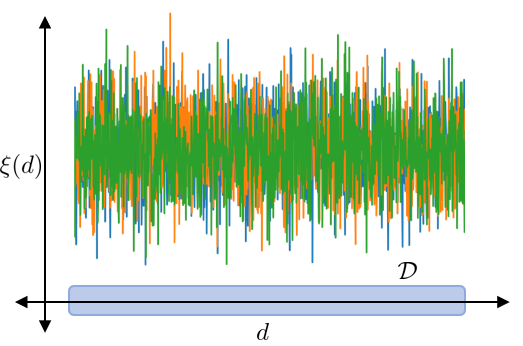}
        \caption{Uncorrelated}
        \label{fig:uncorrelated_field}
    \end{subfigure}
    \caption{Realizations of random field uncertainty $\xi(d)$ obtained from three distinct field types. Notice the static case corresponds to $\xi(d) = \xi$ and can be classified as a random variable (having no dependence on $d$).}
    \label{fig:field_uncertainty}
\end{figure}

Random initial/boundary conditions can be described by the special case where the underlying uncertainty is static $\xi(d) = \xi$ (if needed, implicit dependencies on $\xi$ can be converted into random initial/boundary conditions via augmenting the system of differential equations \cite{soong1973random}):
\begin{equation} 
	\begin{gathered}
		Dy(d, \xi) = f(y(d, \xi), d) \\
		q(y(d_0, \xi), \xi, d_0) = 0, \ d_0 \in \mathcal{D}_0.
	\end{gathered}
	\label{eq:model_random_initial}
\end{equation}
Figure \ref{fig:field_uncertainty} illustrates this distinction for realizations of different random field classifications. These pertain to cases where the underlying uncertainty can be described by a random variable instead of a random field. This often occurs in applications with uncertain parameters that are domain invariant (e.g., kinetic rate constants). 

RDEs that feature inhomogeneous random field input $h(d)\xi(d)$ can be denoted as:
\begin{equation}
	\begin{gathered}
		Dy(d, \xi(d)) = f(y(d, \xi(d)), d) + h(d) \xi(d) \\
		q(y(d_0, \xi(d_0)), d_0) = 0, \ d_0 \in \mathcal{D}_0
	\end{gathered}
	\label{eq:inhomo_rde}
\end{equation}
where $h : \mathcal{D} \mapsto \mathbb{R}$ is a deterministic $\mathbb{R}$-valued function. These typically denote deterministic differential equations $Dy(d) = f(y(d), d)$ that are subjected to random noise $h(d) \xi(d)$. For dynamical systems, these reduce to the well-studied area of SDEs:
\begin{equation}
	\begin{gathered}
		dy(t, \xi(t)) = f(y(t, \xi(t)), t)dt + h(t)d\xi(t) \\
		y(0, \xi(0)) = y_0
	\end{gathered}
	\label{eq:time_sde}
\end{equation}
where $y_0 \in \mathbb{R}^{n_y}$ is the initial value and $d\xi(t) := \xi(t + dt) - \xi(d)$ is the differential of the random field $\xi(t)$ (a stochastic process) which is defined in accordance with an appropriate stochastic integral \cite{graham2013modeling, cinlar2013introduction}. Here, $\xi(t)$ typically denotes a Wiener process (i.e., white noise or Brownian motion), but can be generalized to any process that is a semi-martingale (i.e., it can be represented as a sum of local martingale processes) which includes sample function differentiable, Wiener, and Poisson processes \cite{cinlar2013introduction}. Such SDEs are behind important applications such as the random diffusive and kinetic systems, discussed in Section \ref{sec:intro}. 
\\

The differential equation in \eqref{eq:time_sde} can be ambiguous and is generally understood to denote the integral equation:
\begin{equation}
	y(t + s, \xi(t + s)) - y(t, \xi(t)) = \int_t^{t+s} f(y(t, \xi(t)), t)dt + \int_t^{t+s} h(t)d\xi(t)
	\label{eq:integral_sde}
\end{equation}
where the first integral is a classical Riemann or Lebesgue integral and the second is a stochastic integral (typically an It\^{o} or Stratonovich integral). It\^{o} integrals are the most common choice since they integrate in a non-anticipative manner (making them well-suited for time evolutions) \cite{ito1984foundations}, while Stratonovich integrals are better suited for general topological manifolds (e.g., the input domains of general random fields) \cite{stratonovich1966new}.

The most difficult types of models \eqref{eq:model_rdes} correspond to ones where $\xi(d)$ is used as a random coefficient. These encompass a wide variety of RDEs which compounds the complexity in analysis to establish properties (e.g., the existence and uniqueness of solutions); however, these have been of increased interest in mathematical research \cite{neckel2013random}. An illustrative application would be a transient diffusive process that is subject to spatially random diffusivity $\xi(x)$:
\begin{equation}
	\begin{gathered}
		\frac{\partial y_c(t, x, \xi(x))}{\partial t} = \xi(x) \nabla^2_x y_c(t, x, \xi(x)) + y_g(t, x) \\
		y_c(0, x, \xi(x)) = y_{c,t0} \\
		y_c(t, x_0, \xi(x_0)) = y_{c,x0}, \ x_0 \in \mathcal{D}_{x0}
	\end{gathered}
	\label{eq:transient_diffusion}
\end{equation}
where $y_c(\cdot)$ denotes the concentration/temperature, $y_g(\cdot)$ is a deterministic generative/reactive term, $y_{c,t0} \in \mathbb{R}$ is the initial concentration/temperature field, and $y_{c,x0} \in \mathbb{R}$ are the boundary concentration/temperature values. Such a system is the focus of the case study presented in Section \ref{sec:diffusion_case}.

\subsection{Stochastic Simulation} \label{sec:model_sim}

Now that we can characterize random system models following our discussion in Section \ref{sec:model_characterize}, we can establish how these models can be handled computationally via simulation. The idea behind converting infinite-dimensional models into finite-dimensional forms via transformation follows from the unifying abstraction for InfiniteOpt problems proposed in \cite{pulsipher2021unifying}. Typical transformations can be broken down into three groups: MC, method of weighted residuals (MWR), and probability density  function (pdf) solutions. Our focus in this work will be on MC, since these techniques are frequently used in software tools and are the most intuitive for incorporation into optimization formulations. However, we also provide some discussion with respect to MWR and pdf techniques for completeness.

\subsubsection{Monte Carlo} \label{sec:mc_model}

We begin by considering the simulation of \eqref{eq:model_random_initial} in the special case that the underlying uncertainty can be expressed as a random variable $\xi(d) = \xi$ (i.e., we have a random initial/boundary condition model). With MC samples $\{\hat{\xi}_k: k \in \mathcal{K}\}$, we can express \eqref{eq:model_random_initial} as a system of $|\mathcal{K}|$ deterministic equations:
\begin{equation}
	\begin{gathered}
		D\hat{y}_k(d) = f(\hat{y}_k(d), d), \ k \in \mathcal{K} \\
		q(\hat{y}_k(d_0), \hat{\xi}, d_0) = 0, \ d_0 \in \mathcal{D}_0, k \in \mathcal{K} 
	\end{gathered}
	\label{eq:mc_initial_cond}
\end{equation}
where $\hat{y}_k(d)$ denote the sample functions of the model random field $y(d, \xi)$ that correspond to an MC sample $\xi_k$. These denote a set of standard deterministic differential equations that can be simulated via classical methods (e.g., finite difference or orthogonal collocation over finite elements). We refer the reader to \cite{graham2013modeling, biegler2010nonlinear} for a review of these deterministic solution methods. In practice, these can be simulated numerically with software packages such as \texttt{DifferentialEquations.jl} in \texttt{Julia} or \texttt{odeint} in \texttt{Python}. The ensemble of MC solutions $\{\hat{y}_k(d) : k \in \mathcal{K}\}$ can then be used to approximate the moments of the random field solution $y(d, \xi(d))$ following the approaches detailed in \cite{christakos2017spatiotemporal} as described in Section \ref{sec:field_characteriations}.
\\

A straightforward MC simulation approach for general RDAEs in \eqref{eq:random_daes} involves drawing sample function realizations $\{\hat{\xi}_k(d): k \in \mathcal{K}\}$ from the underlying random field uncertainty $\xi(d)$. With these we can produce $|\mathcal{K}|$ deterministic sets of the DAEs from \eqref{eq:random_daes}:
\begin{equation}
	g(D\hat{y}_k, \hat{y}_k(d), \hat{\xi}_k(d), d) = 0, \ k \in \mathcal{K}.
	\label{eq:mc_random_daes}
\end{equation}
These denote deterministic DAEs that we can solve using traditional techniques in analogous fashion to the deterministic ordinary differential equations \eqref{eq:mc_initial_cond}. One general approach is to model derivatives $D\hat{y}_k$ as auxiliary variables such that \eqref{eq:mc_random_daes} can be treated as an algebraic set of equations in conjunction with a set of ODEs that relate the auxiliary derivative variables to the derivatives $D\hat{y}_k$ \cite{pulsipher2021unifying,biegler2010nonlinear}. Moreover, deterministic DAEs can be readily simulated using software packages such as \texttt{DifferentialEquations.jl} and \texttt{InfiniteOpt.jl}. A key caveat of this MC approach for solving \eqref{eq:random_daes} is that $\xi(d)$ must be sample function continuous and differentiable such that we can avoid using stochastic differentials/integrals that require stochastic calculus techniques and cannot be readily collapsed into a deterministic representation. Moreover, each sample DAE set $g(D\hat{y}_k, \hat{y}_k(d), \hat{\xi}_k(d), d) = 0$ is anticipative over the domain $\mathcal{D}$ which may not be desirable for certain domains (e.g., time). Developing more advanced solution methods (that can avoid the above two caveats) for simulating \eqref{eq:random_daes} is an open-research area, since general RDAE theory is not yet well-developed \cite{suthar2021explicit}. Hence, this MC approach is utilized in most of the case studies presented in Section \ref{sec:cases}.
\\

For system models that use RDEs with inhomogeneous noise terms, such as \eqref{eq:inhomo_rde}, we can employ SDE and SPDE numerical techniques to simulate them. Dynamic SDEs are a popular problem class following the form \eqref{eq:integral_sde} with It\^{o} integrals; these can be approximated with finite difference methods such as Euler–Maruyama, which adapts the implicit Euler method to suit SDEs with Wiener noise $D\xi(t) = dW(t)$ \cite{graham2013modeling}. For a set of discrete times $\{t_i : 1, \dots,N\}$ the $k^\text{th}$, the MC solution is computed as:
\begin{equation}
	\hat{y}_k(t_{i+1}) = \hat{y}_k(t_{i}) + f(\hat{y}_k(t_{i}), t_i)(t_{i+1} - t_i) + h(t_i)(W(t_{i+1})-W(t_i)), \ i = 1,\dots, N-1
	\label{eq:ito_euler}
\end{equation}
where we have that $(W(t_{i+1})-W(t_i)) \sim \mathcal{N}(0, t_{i+1} - t_i)$ and is sampled at each time step. This is the standard approach for simulating SDEs (often called Brownian dynamics simulation) and is non-anticipative (but it is only mildly accurate). High-order (more accurate) methods can also be implemented and these techniques can be adapted for more general SPDEs. We refer the reader to review \cite{zhang2017numerical} for a thorough discussion on the properties and methodologies for implementing these techniques. 

\subsubsection{Other Transformations} \label{sec:other_trans}

\begin{figure}[!htb]
	\includegraphics[width=0.95\textwidth]{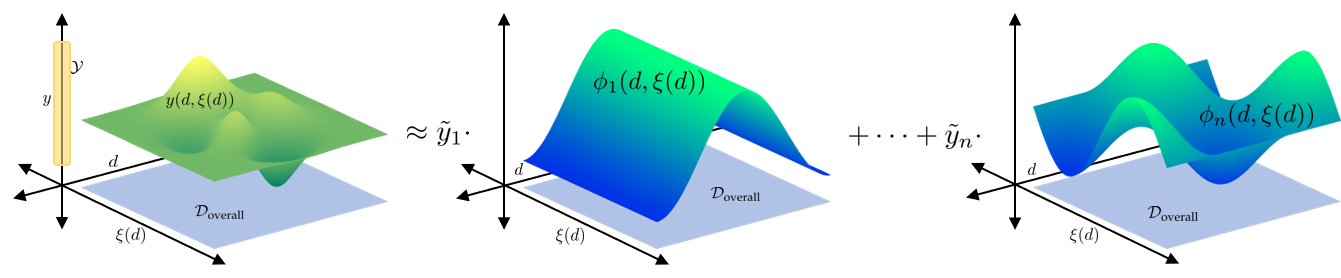}
	\centering
	\caption{Depiction of an MWR approximation for $y(d, \xi(d))$ using a linear combination of orthogonal basis functions $\phi(d, \xi(d))$.}
	\label{fig:mwr_basis}
\end{figure}

Despite the popularity of MC techniques, a number of other transformation techniques have been developed to simulate stochastic systems. The method of weighted residuals denotes a general class of transformations where the system variables are expressed via basis expansions. In the context of stochastic models, MWR transformations are often referred to as polynomial chaos expansion which emphasizes how the functional basis helps to propagate the randomness (chaos) through the model equations. We approximate the system variables with a linear combination of basis functions $\{\phi_i(d, \xi(d)): i \in \mathcal{I}\}$:
\begin{equation}
	y(d, \xi(d)) \approx \sum_{i \in \mathcal{I}} \tilde{y}_i \phi_i(d, \xi(d))
	\label{eq:mwr_basis}
\end{equation}
where $\tilde{y}_i \in \mathbb{R}$ are the basis function coefficients. Typically, the basis functions are chosen such that they are orthogonal relative to the finite distribution functions of $\xi(d)$. Figure \ref{fig:mwr_basis} provides a visual depiction of \eqref{eq:mwr_basis}. Note that this approximation becomes exact as $|\mathcal{I}| \rightarrow \infty$ if $y(d, \xi(d))$ resides in the same function space as the basis functions. We can then project the modeling equations onto the basis functions to yield a system of deterministic algebraic equations with variables $\tilde{y}$ that can be solved using classical methods. We refer interested readers to \cite{zhang2017numerical, youssef2021poly} to learn more about performing MWR transformation techniques to simulate stochastic models.

Another group of transformation methods involve representing the stochastic model in terms of its pdf to obtain a deterministic representation that can be solved using classical approaches. For instance, the Fokker-Planck equation for stochastic particle motion is derived by considering the expectation of an SDE and then applying a number of simplifications to yield a deterministic PDE in terms of the pdf. We refer the reader to \cite{graham2013modeling} to learn more.

\section{Random Field Optimization} \label{sec:optimization}

In this section, we discuss how random field representations along with the stochastic system models that they produce can be incorporated in a general optimization framework. In particular, we discuss how these stochastic InfiniteOpt problems can be formulated and how these formulations can be solved via finite transformations building on the discussion provided in Section \ref{sec:modeling}. 

\subsection{Characterization}

Following work in \cite{pulsipher2021unifying}, we extend the unifying abstraction for InfiniteOpt problems to directly incorporate random field uncertainty which constitutes our framework for random field optimization. Hence, below we outline this framework in this context and establish notation. We invite the reader to see \cite{pulsipher2021unifying} for a rigorous discussion and presentation of this modeling approach.

\subsubsection{Infinite Domain}

\begin{figure}[!htb]
	\includegraphics[width=0.8\textwidth]{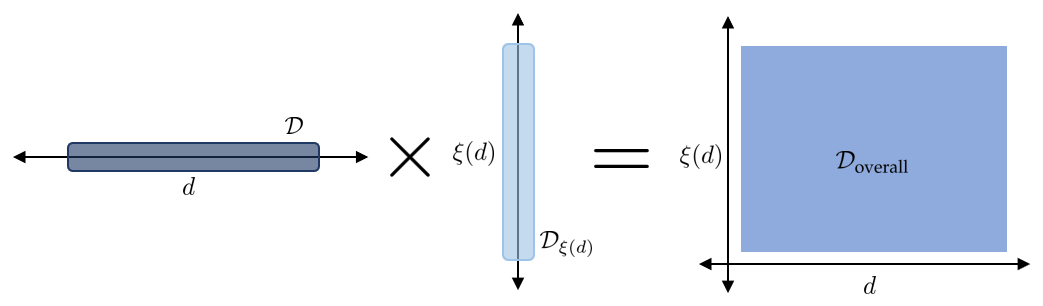}
	\centering
	\caption{An illustration of the overall InfiniteOpt problem domain $\mathcal{D}_{overall}$ defined in Equation \eqref{eq:infiniteopt_domain}.}
	\label{fig:infinite_domain}
\end{figure}

We first consider the infinite domain $\mathcal{D} \subseteq \mathbb{R}^{n_d}$ which denotes the deterministic input domain of the variables (decision functions) of the InfiniteOpt problem. With $\mathcal{D}$, we define the infinite parameter $d \in \mathcal{D}$ which indexes the decision variables. As discussed in Section \ref{sec:random_fields}, we use $d \in \mathcal{D}$ to construct the manifold support for the random field uncertainty $\xi(d) \in \mathcal{D}_{\xi(d)}$ that our system encounters. With this we can then define the overall infinite domain of our InfiniteOpt problem as:
\begin{equation}
	\mathcal{D}_\text{overall} := \mathcal{D} \times \mathcal{D}_{\xi(d)}.
	\label{eq:infiniteopt_domain}
\end{equation}
Figure \ref{fig:infinite_domain} provides a highlevel depiction of $\mathcal{D}_\text{overall}$.

\subsubsection{Variables}

The random field variables $y : \mathcal{D}_\text{overall} \mapsto \mathcal{Y}$ defined in \eqref{eq:model_vars} can be interpreted as a generalization of recourse variables  which are indexed over $\mathcal{D}$. Moreover, we define deterministic decision variables $z : \mathcal{D} \mapsto \mathcal{Z} \subseteq \mathbb{R}^{n_z}$:
\begin{equation}
	z(d) \in \mathcal{Z}.
	\label{eq:opt_vars}
\end{equation}
These typically correspond to parameters and/or deterministic functions in the stochastic models considered in Section \ref{sec:modeling} which constitute the primary decision functions we wish to make. In this sense, the variables $z(d)$ can be interpreted as an infinite-dimensional generalization of first-stage variables associated with classical two-stage stochastic programming formulations. These also capture finite variables $z(d) = z$ as a special case. 

\subsubsection{Objective}

\begin{figure}[!htb]
	\includegraphics[width=0.7\textwidth]{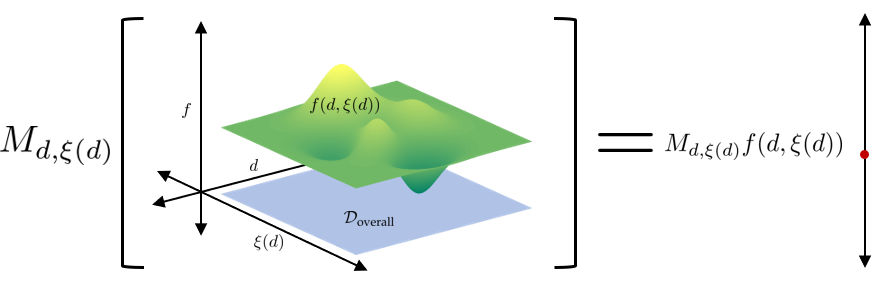}
	\centering
	\caption{Depiction of how an RFO objective cost function $f(d, \xi(d))$ is scalarized via the general measure operator $M_{d, \xi(d)}$.}
	\label{fig:objective}
\end{figure}

Objective functions employ a measure operator $M_{d, \xi(d)} : (\mathcal{D}, \mathcal{D}_{\xi(d)}) \mapsto \mathbb{R}$ which summarizes (scalarizes) an infinite-dimensional cost function $f(Dy, y(d, \xi(d)), z(d), \xi(d), d)$:
\begin{equation}
	\min M f(Dy, y(d, \xi(d)), z(d), \xi(d), d).
	\label{eq:objective}
\end{equation}
For convenience in notation, we will sometimes use $f(d, \xi(d)) := f(Dy, y(d, \xi(d)), z(d), \xi(d), d)$. Figure \ref{fig:objective} illustrates Equation \eqref{eq:objective}. For example, we can use a measure operator $M_{d, \xi(d)}$ that uses the random field expectation $\mathbb{E}_{\xi(d)}$ and Riemann integrals to obtain:
\begin{equation}
	\min \mathbb{E}_{\xi(d)}\left[ \int_{d \in \mathcal{D}} f(Dy, y(d, \xi(d)), z(d), \xi(d), d) d(d) \right].
	\label{eq:objective_example}
\end{equation}
For simplicity in presentation, this serves as the workhorse for the numerical studies presented in Section \ref{sec:cases}. 

However, we can envision developing other measure operators to shape cost functions $f(d, \xi(d))$. For instance, we can employ the excursion probability in \eqref{eq:excursion_prob}:
\begin{equation}
	\min \mathbb{P}_{\xi(d)}\left(\max_{d \in \mathcal{D}} f(d, \xi(d)) \geq u \right)
\end{equation}
which minimizes the probability that the cost function $f(d, \xi(d))$ exceeds $u \in \mathbb{R}$. Another interesting class of measures results from extending classical risk measures from stochastic optimization (SO) to operate on random field cost functions $f(d, \xi(d))$. The conditional-value-at-risk (CVaR) measure operator is popular in SO that seeks to measure the tail of the pdf of the cost function:
\begin{equation}
	\text{CVaR}_\xi(f(\xi); \alpha) := \min_{\hat{f} \in \mathbb{R}}\left\{\hat{f} + \frac{1}{1-\alpha}\mathbb{E}_\xi[\max(f(\xi) - \hat{f}, 0)]\right\}
\end{equation}
where $\xi \in \mathcal{D}_\xi$ is a random variable, $\hat{f} \in \mathbb{R}$ is an auxiliary variable, and $\alpha \in [0, 1)$ is the probability level that defines the portion of the tail that is measured. One way to generalize CVaR to random field cost functions would be to nest a measure $M_d$:
\begin{equation}
	\text{CVaR}_{\xi(d)}(f(d, \xi(d)); \alpha) := \min_{\hat{f} \in \mathbb{R}}\left\{\hat{f} + \frac{1}{1-\alpha}\mathbb{E}_{\xi(d)}[\max(M_df(d, \xi(d)) - \hat{f}, 0)]\right\}
\end{equation}
where $M_df(d, \xi(d))$ can be a traditional deterministic measure such as the integral $\int_{d \in \mathcal{D}} f(d, \xi(d)) d(d)$ or the maximization $\max_{d \in \mathcal{D}} f(d, \xi(d))$. We leave a rigorous exploration of such extensions to future work.

\subsubsection{Constraints}

\begin{figure}[!htb]
	\includegraphics[width=0.4\textwidth]{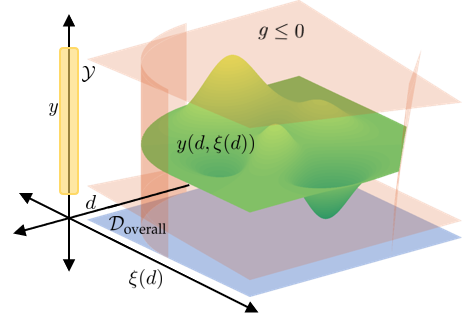}
	\centering
	\caption{Depiction of an algebraically-constrained RFO problem following \eqref{eq:dae_constraints}.}
	\label{fig:constraints}
\end{figure}

We establish a couple of constraint classes, algebraic and measure constraints, which encode modeling equations (i.e., those presented in Section \ref{sec:modeling}) and decision restrictions. Algebraic constraints are established by generalizing the random DAE presented in \eqref{eq:random_daes}:
\begin{equation}
	g(Dy, y(d, \xi(d)), z(d), \xi(d), d) \leq 0, \; \xi(d) \in \mathcal{D}_{\xi(d)}, d \in \mathcal{D}.
	\label{eq:dae_constraints}
\end{equation}
This enforces the constraint $g \leq 0$ over each realization of $(\xi(d), d)$. Moreover, this captures point constraints (e.g., initial/boundary conditions) by enforcing this constraint and restricting the overall infinite domain $\mathcal{D}_\text{overall}$ to a particular point $(\xi(d^*), d^*)$. We can also capture the more structured R(P)DE equation forms presented in Section \eqref{sec:model_characterize} as special cases. We provide a high-level representation of a simple algebraic constraint system in Figure \ref{fig:constraints}.

Constraints that we enforce almost surely with respect to the random field uncertainty $\xi(d)$ can be overly restrictive in certain applications, producing excessively conservative solutions. Measure constraints help alleviate this limitation by enforcing a constraint on a measure function:
\begin{equation}
	Mh(y(d, \xi(d)), z(d), \xi(d), d) \leq 0
	\label{eq:measure_constraints}
\end{equation}
where $M$ is an arbitrary measure operator and $h(\cdot)$ is an $\mathbb{R}^{n_h}$-valued function. This constrains a summarizing statistic over the random field model rather than enforcing a condition on every function sample realization. For instance, we can adapt the excursion probability in \eqref{eq:excursion_prob} to constrain the probability of a variable $y(d, \xi(d))$ exceeding a certain threshold $u$ to a probability level $\alpha \in [0, 1]$:
\begin{equation}
	\mathbb{P}_{\xi(d)}\left(\max_{d \in \mathcal{D}} y(d, \xi(d)) \geq u \right) \leq \alpha.
	\label{eq:excursion_constraint}
\end{equation}
This excursion probability constraint is a generalization of classical chance constraints uses in SO. Constraint \eqref{eq:excursion_constraint} can be reformulated:
\begin{equation}
	\mathbb{E}_{\xi(d)}\left[\mathbbm{1}(\max_{d \in \mathcal{D}} y(d, \xi(d)) \geq u)\right] - \alpha \leq 0
\end{equation}
which conforms to \eqref{eq:measure_constraints} by setting $My = \mathbb{E}[\mathbbm{1}(\max_d y \geq u)] - \alpha$. We note that when $\alpha \rightarrow 0$ we seek to enforce the constraint:
\begin{equation}
	\max_{d \in \mathcal{D}} y(d, \xi(d)) < u, \ \xi(d) \in \mathcal{D}_{\xi(d)}
\end{equation}
which is equivalent to:
\begin{equation}
	y(d, \xi(d)) < u, \ \xi(d) \in \mathcal{D}_{\xi(d)}, d \in \mathcal{D}.
\end{equation}
Hence, the special case of \eqref{eq:excursion_constraint} corresponds to an algebraic constraint in accordance with \eqref{eq:dae_constraints}. Section \ref{sec:diffusion_case} implements this proposed constraint class in an illustrative study.

\subsubsection{RFO Formulation}

We incorporate all the above modeling objects to define the general RFO problem:
\begin{equation}
	\begin{aligned}
		&&&\min && Mf(Dy, y(d, \xi(d)), z(d), \xi(d), d) \\
		&&& \text{s.t.} && g(Dy, y(d, \xi(d)), z(d), \xi(d), d) \leq 0, && \xi(d) \in \mathcal{D}_{\xi(d)}, d \in \mathcal{D} \\
		&&&&& Mh(y(d, \xi(d)), z(d), \xi(d), d) \leq 0.
	\end{aligned}
	\label{eq:general_form}
\end{equation}
This summarizes modeling elements that an RFO problem may entail; moreover, we highlight that it is a special case of the general InfiniteOpt problem formulation presented in \cite{pulsipher2021unifying}.

\subsection{Transformations}

Random field optimization problems are not readily amenable to classical optimization solution methods due to their infinite-dimensional nature. Hence, we seek to obtain approximate finite representations of RFO problems via transformation routines.  

\subsubsection{Direct Transcription} \label{sec:direct_transcription}

We can use MC sampling approaches to project the RFO problem equations onto a deterministic basis that can be solved with classical methods (e.g., direct transcription). The most straightforward approach is to collect MC sample functions $\hat{\mathcal{D}}_{\xi(d)} := \{\hat{\xi}_k(d) : k \in \mathcal{K}\}$ and define a set of sample points $\hat{\mathcal{D}} := \{\hat{d}_i : i \in \mathcal{I}\}$. Projecting the variables $y(d, \xi(d))$ and $z(d)$ onto these sets provides a collection of finite variables:
\begin{equation}
	\begin{aligned}
		&y(\hat{d}, \hat{\xi}_k(\hat{d})) \in \mathcal{Y}, && \hat{\xi(d)} \in \hat{\mathcal{D}}_{\xi(d)}, \hat{d} \in \hat{\mathcal{D}} \\
		&z(\hat{d}) \in \mathcal{Z}, && \hat{d} \in \hat{\mathcal{D}} 
	\end{aligned}
	\label{eq:finite_vars}
\end{equation}
which we can equivalently write as $\{y_{ik} \in \mathcal{Y}: i \in \mathcal{I}, k \in \mathcal{K}\}$ and $\{z_{i} \in \mathcal{Z}: i \in \mathcal{I}\}$. This projection onto discrete points is called direct transcription.
\\

We can use these variables to approximate the modeling equations and constraints using the methods discussed in Section \ref{sec:mc_model}. For instance, we can use the Euler–Maruyama method if the equations are RDEs with a non-homogeneous noise term to obtain a finite set of equations. Algebraic equations/constraints can be projected directly onto $\hat{\mathcal{D}}_{\xi(d)}$ and $\hat{\mathcal{D}}$:
\begin{equation}
	g(y_{ik}, z_i, \hat{\xi}_{ik}, \hat{d}_i) \leq 0, \; i \in \mathcal{I}, k \in \mathcal{K}.
\end{equation}
Similarly, we can project random DAE constraints directly by introducing auxiliary variables $Dy_{ik}$ for the derivatives $Dy(\hat{d}_i, \hat{\xi}_k(\hat{d}_i))$ which can be approximated using the support points with deterministic methods (e.g., explicit Euler) \cite{pulsipher2021unifying}. However, recall that this approach is anticipative and requires that the random field $\xi(d)$ be sample function differentiable. 
\\

Measures can be approximated directly using the sample functions/points. For instance, the objective function in \eqref{eq:objective_example} would be expressed:
\begin{equation}
	|\mathcal{K}|^{-1}\sum_{k \in \mathcal{K}} \sum_{i \in \mathcal{I}} \beta_i f(Dy_{ik}, y_{ik}, z_i, \hat{\xi}_k(\hat{d}_i), \hat{d}_i)
\end{equation}
where $\beta_i, \ i \in \mathcal{I}$ are constants (e.g., quadrature coefficients) chosen to approximate the integral measure. For transforming excursion probability measures, the analytical approaches described in Section \ref{sec:fields_measures} (e.g., \eqref{eq:rice_formula}) for computing the excursion probability may not always be useful due to the difficulty in propagating the moments of random fields through an optimization formulation. In some cases,  \eqref{eq:mc_expected_ec} can be used if the explicit form of $\mathcal{X}(A_u(\cdot)))$ is amenable to the desired optimization tools (e.g., if it is compatible with auto-differentiation to enable gradient optimization solution methods). Otherwise, we can approximate the excursion probability via our MC samples:
\begin{equation}
	P_{\xi(d)}\left\{\sup_{d \in \mathcal{D}} f(d, \xi(d)) \geq u\right\} \approx |\mathcal{K}|^{-1} \sum_{k \in \mathcal{K}} \mathbbm{1}\left(\sup_{i \in \mathcal{I}} f(\hat{d}_i,\hat{\xi}_k(\hat{d}_i)) \geq u\right).
	\label{eq:mc_excursion_prob2}
\end{equation}
Here, the indicator function $\mathbbm{1}(\cdot)$ can be treated as a nonlinear function, or we can reformulate \eqref{eq:mc_excursion_prob2} into $|\mathcal{K}|^{-1} \sum_{k \in \mathcal{K}} q_{k}$ using big-M constraints:
\begin{equation}
	\begin{aligned}
	&\overline{f}_k \geq f(\hat{d}_i,\hat{\xi}_k(\hat{d}_i)), && i \in \mathcal{I}, k \in \mathcal{K} \\
	&\overline{f}_k - u\leq q_{k}M, && k \in \mathcal{K} 
	\end{aligned}
	\label{eq:big_m_escursion_prob}
\end{equation}
where $q_k \in \{0,1\}$ are binary variables and $\overline{f}_k \in \mathbb{R}$ are auxiliary variables. Following the relaxation approach proposed in \cite{pulsipher2019scalable}, we can relax the binary variables to $q_k \in [0, 1]$ (yielding a more tractable continuous formulation) and apply a rounding rule with the optimal values $q_k^*$ which often recovers the binary solution to high accuracy. 
\\

We note that transforming an RFO problem via direct transcription can quickly become intractable for even a moderate number of sample functions $\hat{\xi}_k$ and points $\hat{d}_i$ due to the inherent combinatorics involved in constructing the discretization grid. Decomposition approaches can help to alleviate this limitation such as the overlapping Schwarz decomposition method proposed in \cite{shin2020decentralized}.

\subsubsection{Alternative Transformations} \label{sec:rfo_other_transforms}

Drawing inspiration from Section \ref{sec:other_trans}, we can envision using more advanced transformation methods to obtain a finite problem formulation that can be used with classical optimization solution techniques. For instance, we can use generalized chaos expansion (i.e., MWR) approaches to project the problem over basis functions $\{\phi_i(d, \xi(d)): i \in \mathcal{I}\}$ as shown in Equation \eqref{eq:mwr_basis}. We can then use Galerkin projection to represent Problem \eqref{eq:general_form} in terms of finite coefficient variables to obtain a finite formulation. A similar approach is employed in \cite{chen2010level} to project an RFO problem for topological design into a finite formulation. Developing a general MWR transformation framework for RFO problems is beyond the scope of this work, but denotes a promising area of future research that can build upon the general MWR approach for InfiniteOpt problems presented in \cite{pulsipher2021unifying}.
\\

As noted in Section \ref{sec:modeling}, developing theory and solution methods for general random DAEs is an active area of research. Further research in this area will be invaluable for developing more advanced transformation techniques to apply RFO problems.

\section{Case Studies} \label{sec:cases}

We now present case studies that showcase the applicability of the proposed modeling abstraction. These further motivate the advantage of capturing uncertainty over general domains using our proposed framework. All of these problems are solved using MC approaches implemented using \texttt{InfiniteOpt.jl v0.5.1} and \texttt{GaussianRandomFields.jl v2.1.4} \cite{pulsipher2021unifying}. All the scripts and data needed for reproducing the results are available at \url{https://github.com/zavalab/JuliaBox/tree/master/RandomFieldOptCases}.

\subsection{Simple Dynamical System} \label{sec:power_network}

\begin{figure}[!htb]
	\includegraphics[width=0.5\textwidth]{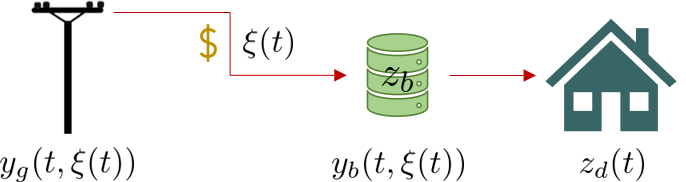}
	\centering
	\caption{Schematic of the power network of case study in Section \ref{sec:power_network}.}
	\label{fig:power_network}
\end{figure}

We use a dynamic optimal design study to illustrate a simple RFO problem that involves dynamic uncertainty. In particular, we consider a home that participates in a dynamic electricity pricing market with random field energy pricing $\xi(t) \in \mathbb{R}_+$, where we seek to optimally design the maximum capacity $z_b \in \mathbb{Z}_b \subseteq \mathbb{R}_+$ of a storage battery to fully take advantage of the market such that we minimize our daily electricity cost. Figure \ref{fig:power_network} details the linear network where electrical power $y_g(t, \xi(t)) \in \mathbb{R}_+$ is purchased from the grid with a unit price of $\xi(t)$ to charge a battery with capacity $y_b(t, \xi(t)) \in [0, z_b]$ that powers a home with electricity demand $z_d(t) \in \mathbb{R}_+$.

\begin{figure}[!htb]
	\includegraphics[width=0.4\textwidth]{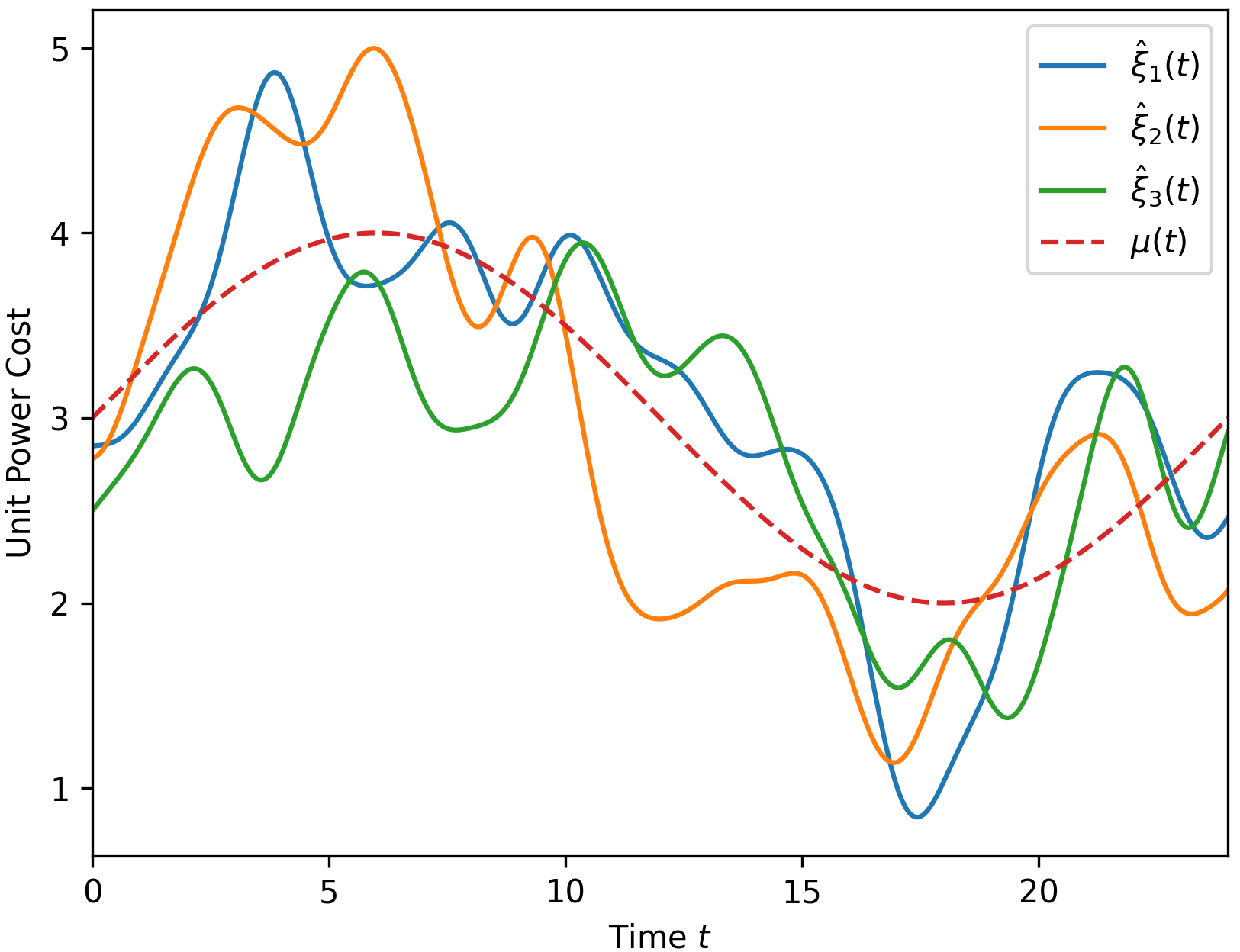}
	\centering
	\caption{Sample function realizations $\hat{\xi}_k(t)$ for $\xi(t)$ following the moments given in \eqref{eq:power_moments}.}
	\label{fig:power_fields}
\end{figure}

The random field cost function $\xi(t)$ for a 24-hour period is characterized as a Gaussian random field with moments:
\begin{equation}
	\begin{aligned}
		&\mu(t) = \sin\left(\frac{2\pi t}{24}\right) + 3 \\
		&\Sigma_\text{SE}(t, t') = \sigma\exp\left(-\frac{(t-t')^2}{2\beta^2}\right)
	\end{aligned}
	\label{eq:power_moments}
\end{equation}
where we set $\sigma = 0.6$ and $\beta=1.5$. The field is modeled with \texttt{GaussianRandomFields.jl} and Figure \ref{fig:power_fields} shows three sample functions $\hat{\xi}_k(t)$ and the mean function $\mu(t)$.

\begin{figure}[!htb]
	\includegraphics[width=0.4\textwidth]{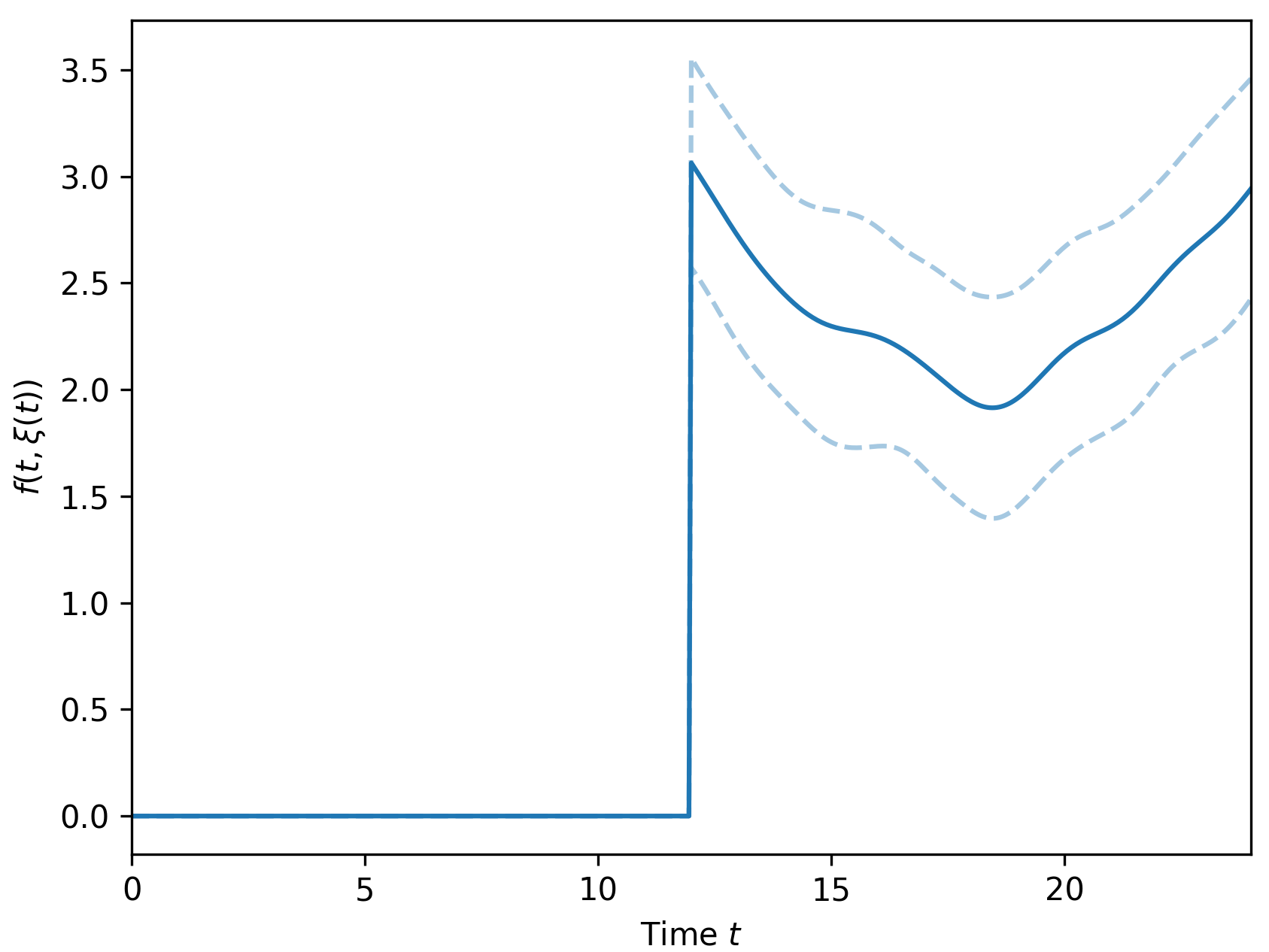}
	\centering
	\caption{Simulated MC ensemble of the operational cost $f(t, \xi(t))$.}
	\label{fig:power_sim}
\end{figure}

The expected total cost incurred is computed via operating on the cost function $f(t, \xi(t)) = \xi(t) y_g(t, \xi(t)) + 0.3z_b$ with $M_{t,\xi(t)}$ from Equation \eqref{eq:objective_example}:
\begin{equation}
	\mathbb{E}_{\xi(t)}\left[\int_{t \in \mathcal{D}_t} \xi(t) y_g(t, \xi(t))dt\right] + 0.3z_b
	\label{eq:power_objective}
\end{equation}
where $\mathcal{D}_t = [0, 24]$. We model the power network by imposing a power balance at the battery node:
\begin{equation}
	\frac{\partial y_b(t, \xi(t))}{\partial t} = y_g(t, \xi(t)) - z_d(t), \ \xi(t) \in \mathcal{D}_{\xi(t)}, t \in \mathcal{D}_t.
	\label{eq:power_balance}
\end{equation}
For simplicity in example, we set $z_d(t) = 1$. Equation \eqref{eq:power_balance} is coupled with periodic boundary conditions that exact the battery be at 50\% capacity at the beginning and end of a 24-hour period:
\begin{equation}
	\begin{aligned}
		&y_b(0, \xi(0)) = 0.5z_b,  && \xi(t) \in \mathcal{D}_{\xi(t)} \\
		&y_b(24, \xi(24)) = 0.5z_b,  && \xi(t) \in \mathcal{D}_{\xi(t)}.
	\end{aligned}
	\label{eq:power_boundaries}
\end{equation}
We combine \eqref{eq:power_objective}-\eqref{eq:power_boundaries} to simulate the operation of our model for fixed $y_g(t, \xi(t)$ and $z_b$ to derive an ensemble of battery response trajectories and costs. Such a simulation is trivial in this case since \eqref{eq:power_balance} does not explicitly depend on $\xi(t)$ and thus can be treated deterministically for each $y_g(t, \hat{\xi}_k(t)$. Figure \ref{fig:power_sim} shows the simulated ensemble of the cost function $f(t, \xi(t))$ with 
\begin{equation}
	y_g(t, \xi(t) = \begin{cases}0 & t \leq 12 \\ 1 & t > 12 \end{cases}
\end{equation}
that is generated via \texttt{InfiniteOpt.jl} using 100 MC samples of $\xi(t)$, 49 equidistant time points, and an implicit Euler method.

We proceed to characterize the RFO problem to optimally choose $y_g(t, \xi(t)$ and $z_b$. We enforce that the battery capacity be kept above 20\% and below 100\% to prevent battery degradation:
\begin{equation}
	0.2z_b \leq y_b(t, \xi(t)) \leq z_b, \ \xi(t) \in \mathcal{D}_{\xi(t)}, t \in \mathcal{D}_t.
	\label{eq:power_restricts}
\end{equation}
With this last constraint, we can combine \eqref{eq:power_objective}-\eqref{eq:power_restricts} to yield the dynamic RFO formulation:
\begin{equation}
	\begin{aligned}
		&&\min_{y_g(\cdot), y_b(\cdot), z_b} &&& \mathbb{E}_{\xi(t)}\left[\int_{t \in \mathcal{D}_t} \xi(t) y_g(t, \xi(t))dt\right] + 0.3z_b\\
		&& \text{s.t.} &&& \frac{\partial y_b(t, \xi(t))}{\partial t} = y_g(t, \xi(t)) - 1, && \xi(t) \in \mathcal{D}_{\xi(t)}, t \in \mathcal{D}_t\\
		&&&&& y_b(0, \xi(0)) = 0.5z_b,  && \xi(t) \in \mathcal{D}_{\xi(t)} \\
		&&&&& y_b(24, \xi(24)) = 0.5z_b,  && \xi(t) \in \mathcal{D}_{\xi(t)} \\
		&&&&& 0.2z_b \leq y_b(t, \xi(t)) \leq z_b, && \xi(t) \in \mathcal{D}_{\xi(t)}, t \in \mathcal{D}_t \\
		&&&&& z_b \in \mathcal{Z}_b
	\end{aligned}
	\label{eq:power_form}
\end{equation}
where we let $\mathcal{Z}_b = [0, 100]$. We implement \eqref{eq:power_form} in \texttt{InfiniteOpt.jl} where we obtain a finite formulation via direct transcription using 49 time points and 100 random field samples $\hat{\xi}_k(t)$ in combination with \texttt{Gurobi v.9.1.2}. Moreover, we compare this solution with the one obtained by solving \eqref{eq:power_form} with deterministic pricing $\mu(t)$. We simulate the system response when using $z_{b,\text{deterministic}}$ against the same random field samples to compute the expected total cost in the deterministic case.

\begin{table}[!htb]
    \centering
    \caption{Optimal solutions for Problem \eqref{eq:power_form} and its deterministic variant.}
    \begin{tabular}{|c | c | c|} 
        \hline
        Formulation &$z_b^*$ & $M_{t, \xi(t)}f^*(t, \xi(t))$ \\
        \hline \hline
		Stochastic & 45.0 & 49.5 \\
		Deterministic & 78.3 & 55.2\\
        \hline
    \end{tabular}
    \label{tab:power_results}
\end{table}

\begin{figure}[!htb]
    \centering
	\begin{subfigure}[b]{0.45\textwidth}
        \centering
        \includegraphics[width=\textwidth]{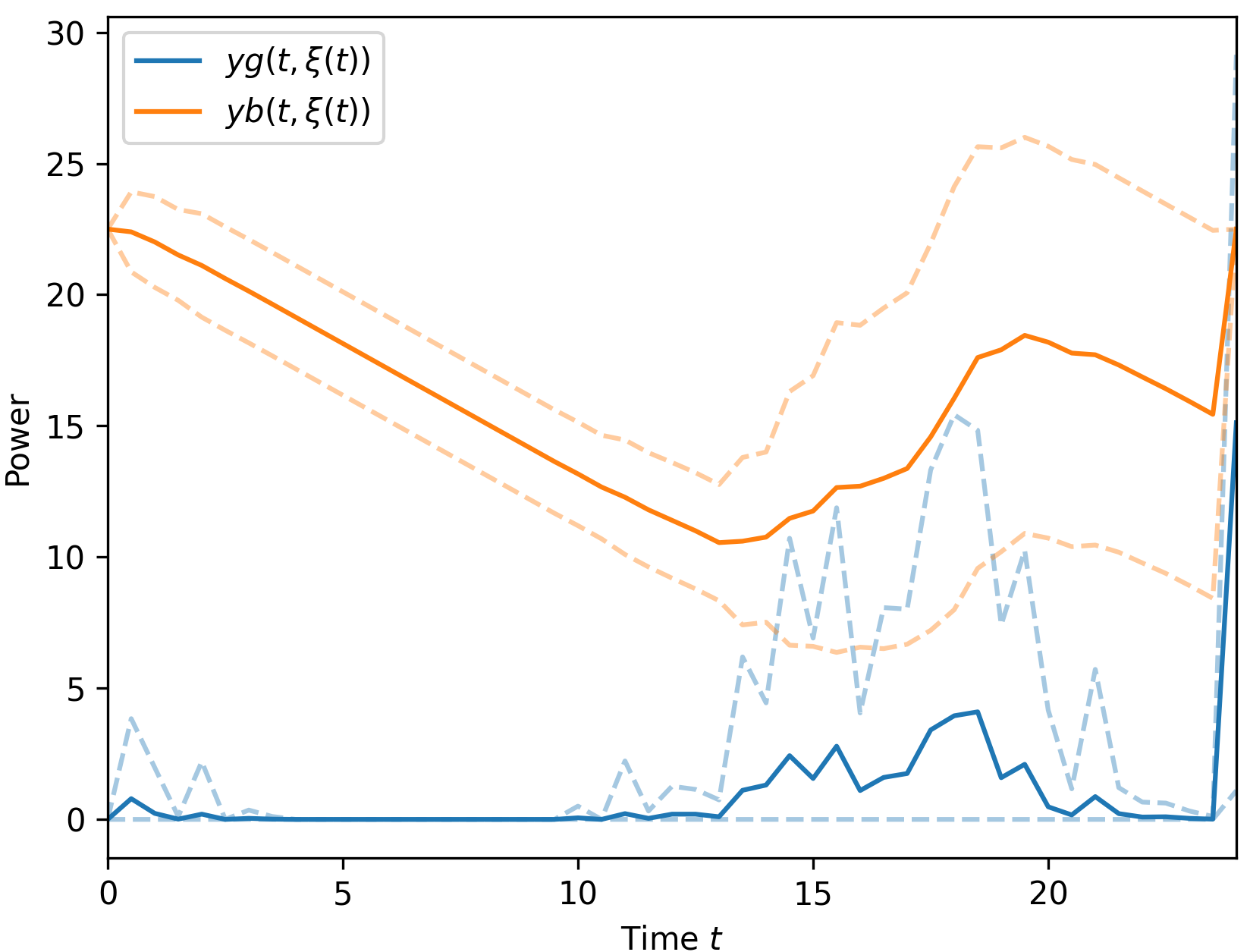}
        \caption{Stochastic}
        \label{fig:stochastic_power}
    \end{subfigure}
	\quad
    \begin{subfigure}[b]{0.45\textwidth}
        \centering
        \includegraphics[width=\textwidth]{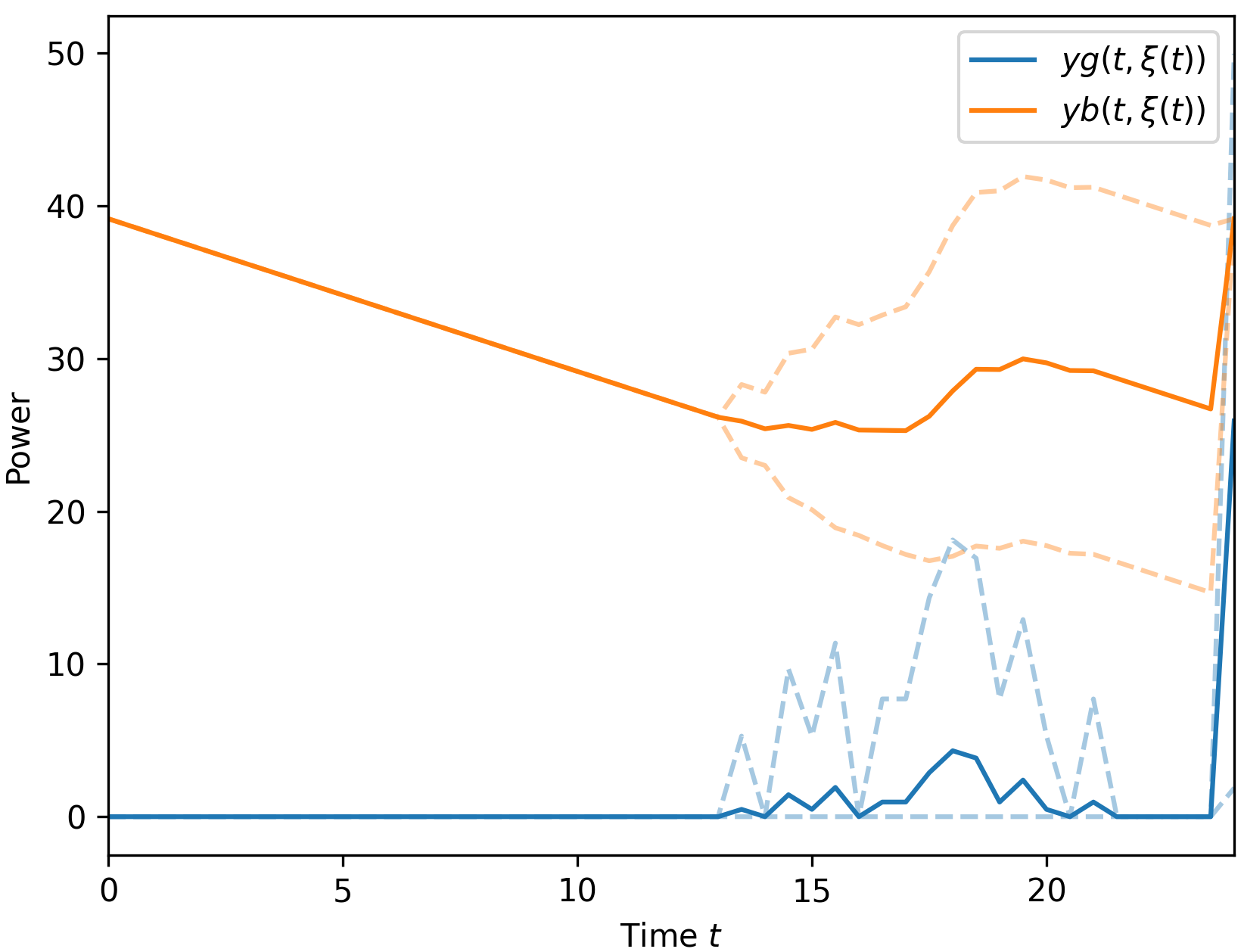}
        \caption{Deterministic}
        \label{fig:deterministic_power}
    \end{subfigure}
    \caption{The ensemble responses for operating the power network over a 24-hour period using the optimal maximum battery capacities $z_b^*$ chosen via the stochastic and deterministic variants of Problem \eqref{eq:power_form}.}
    \label{fig:power_results}
\end{figure}

Table \ref{tab:power_results} shows the results from these studies. The stochastic solution incurs an expected cost that is 10.3\% lower than the deterministic counterpart. This difference can be interpreted as an analogue to the value of stochastic solution (VSS) metric that is commonly used in stochastic programming for evaluating the effectiveness of stochastic solutions against their deterministic counterparts \cite{birge2011introduction}. The cost savings can readily be attributed to how the stochastic formulation selects a more conservative maximum battery capacity $z_b$ that is adequate for the system. Figure \ref{fig:power_results} compares the simulated response of the power network using these choices of $z_b$. We observe how the RFO-derived simulation better utilizes purchasing electricity when it is cheaper later in the period; whereas the deterministic derived simulation uses its larger battery capacity to avoid purchasing to a greater extent.

\subsection{Atomic Layer Deposition under Spatial Uncertainty} \label{sec:ald}

\begin{figure}[!htb]
	\includegraphics[width=0.55\textwidth]{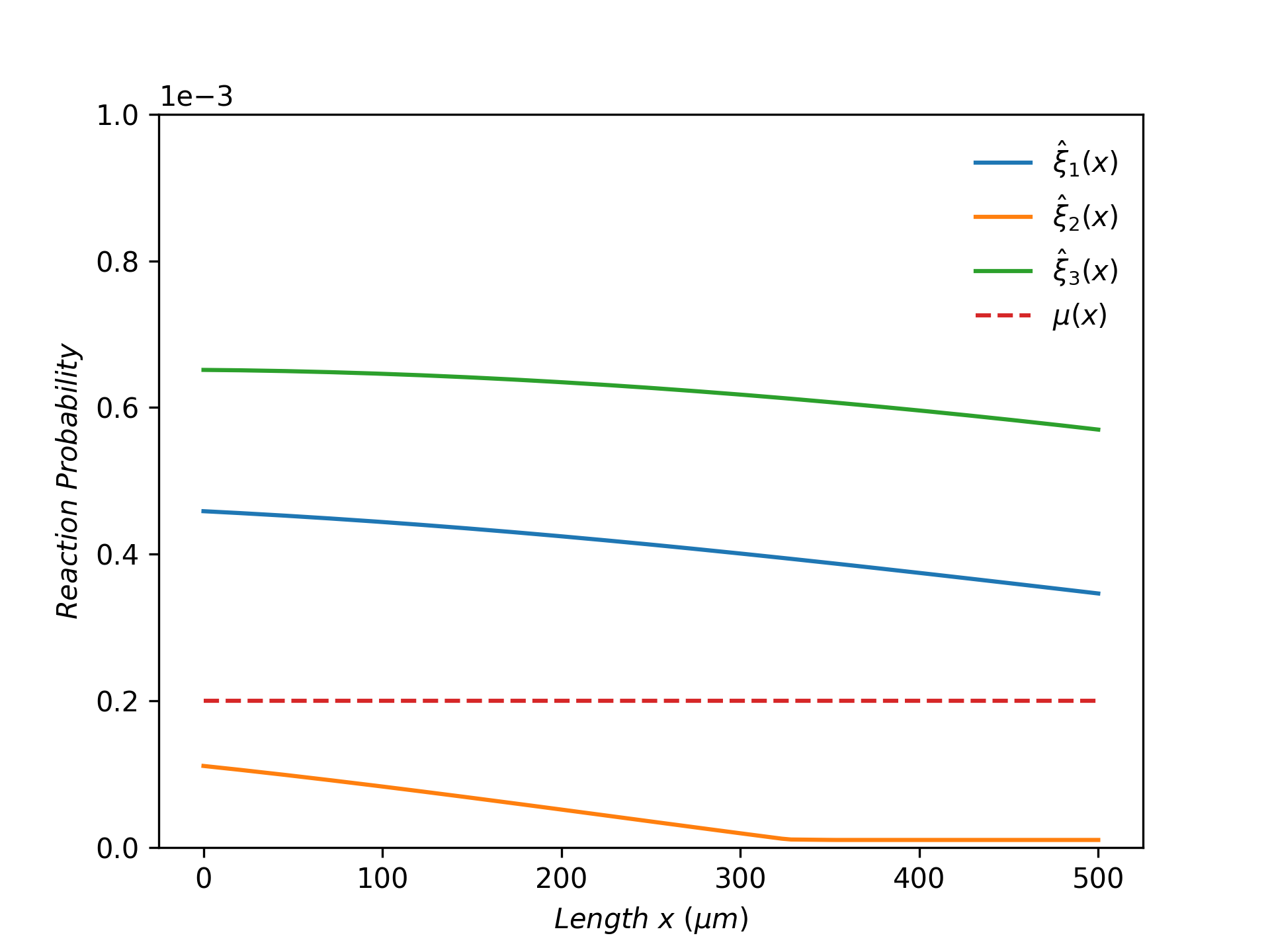}
	\centering
	\caption{Sample function realizations $\hat{\xi}_k(x)$ for uncertain reaction probability of ALD substrate.}
	\label{fig:ald_random_field}
\end{figure}

We now exemplify the application of RFO in problems involving spatial uncertainty. In particular, we examine the atomic layer deposition (ALD) process on a porous substrate surface that is presented in \cite{keuter2015modeling}. This amounts to a 1D reactive diffusion system where a substrate is exposed to a gaseous precursor that diffuses and reacts with it (this process is taken to be isotropic across the surface). We add the complexity of an uncertain reaction probability $\xi(x) : \mathcal{D}_{\xi(x)} \mapsto [0,1]$ that varies along the depth of the reaction substrate; this is taken to be a Gaussian random field with:
\begin{equation}
    \begin{aligned}
        &\mu(x)= 2\times10^{-4} \\
        &\Sigma_{\mathrm{SE}}\left(x, x'\right) = 0.0003 \exp \left(-\frac{\left(x-x'\right)^{2}}{1.8 \times 10^5}\right)
    \end{aligned}
\end{equation}
where we taken $x$ to be in units of $\mu m$. Figure \ref{fig:ald_random_field} depicts realizations of $\xi(x)$ along with its mean function $\mu(x)$.

\begin{figure}[!htb]
	\includegraphics[width=0.55\textwidth]{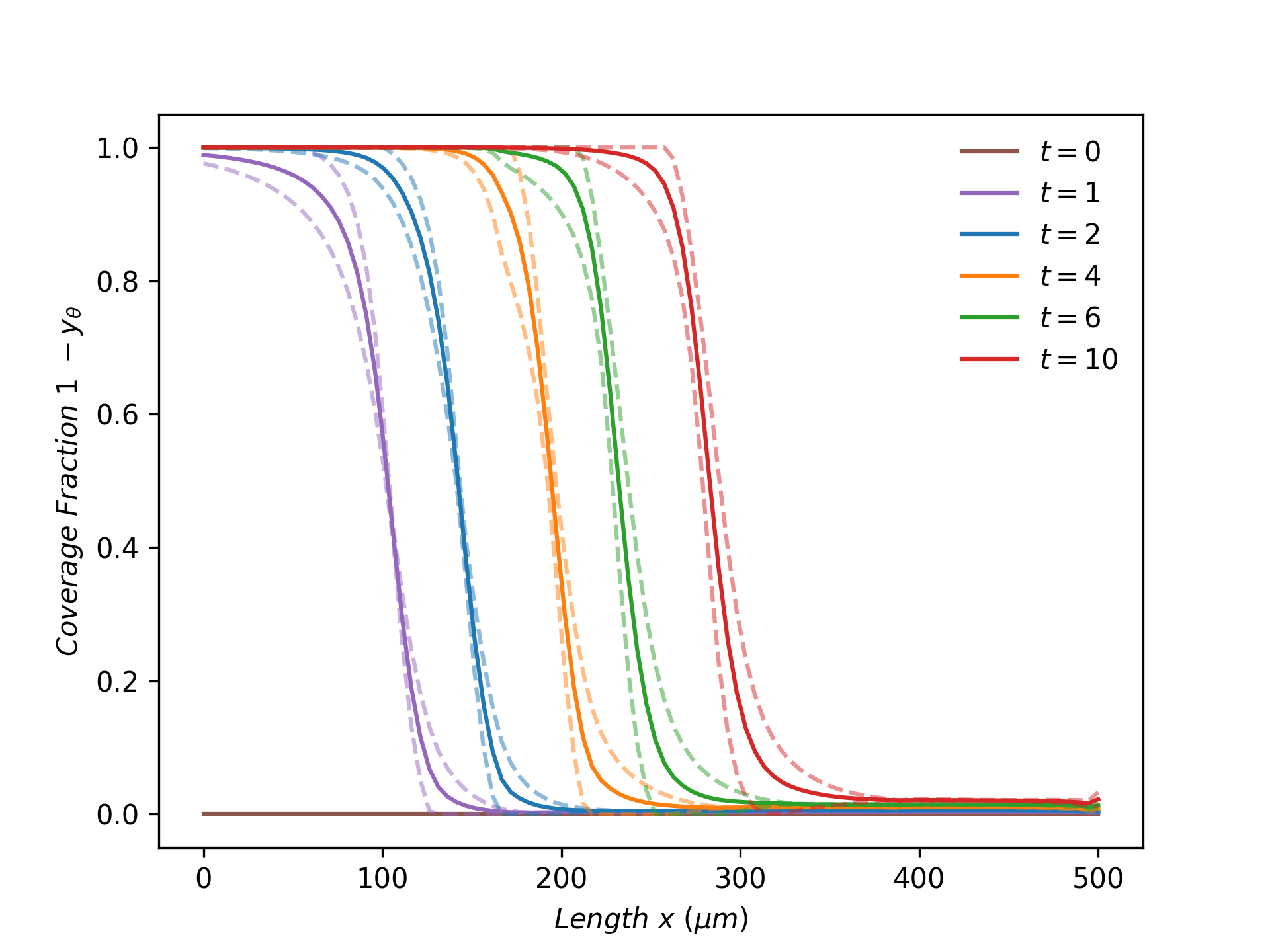}
	\centering
	\caption{Simulated MC ensemble of the coverage profile $1-y_\theta(t, x, \xi(x))$ at select time instances following \eqref{eq:ald_pdes} and \eqref{eq:ald_boundary_conds}.}
	\label{fig:ald_sim}
\end{figure}

The ALD model considers the precursor density $y_p(t, x, \xi(x)) \in \mathbb{R}_+$ in relation to the fraction of available sites in the substrate surface layer $y_\theta(t, x, \xi(x)) \in [0, 1]$. Hence, the surface coverage is given by $1-y_\theta(t, x, \xi(x))$. Following the derivation presented in \cite{keuter2015modeling}, we obtain the the RPDE modeling equations:
\begin{equation}
    \begin{aligned}
        &\frac{\partial y_p(t, x, \xi(x))}{\partial t} = D \frac{\partial^{2} y_p(t, x, \xi(x))}{\partial x^2} - \gamma \xi(x) y_p(t, x, \xi(x)) y_\theta(t, x, \xi(x)), && \xi(x) \in \mathcal{D}_{\xi(x)}, (t,x) \in \mathcal{D}_{t,x} \\
        &\frac{\partial y_\theta(t, x, \xi(x))}{\partial t} = -\eta \xi(x) y_p(t, x, \xi(x)) y_\theta(t, x, \xi(x)), && \xi(x) \in \mathcal{D}_{\xi(x)}, (t,x) \in \mathcal{D}_{t,x}
    \end{aligned}
    \label{eq:ald_pdes}
\end{equation}
where $\mathcal{D}_{t,x} \in \mathcal{D}_t \times \mathcal{D}_x = [0, 10] \times [0, 500]$, $D = 2.81 \times 10^6 \mu m^2s^{-1}$ is the diffusivity constant, $\gamma = 6.912 \times 10^8 s^{-1}$ is the precursor species reaction constant, and $\eta = 1.538 \times 10^7 \mu m^3 s^{-1}$ is the site reaction constant. For boundary conditions we enforce:
\begin{equation}
    \begin{aligned}
        &y_p(0, x,\xi(x)) = 0, && \xi(x) \in \mathcal{D}_{\xi(x)}, x \in \mathcal{D}_{x} \\
        &y_p(t, 0,\xi(0)) = z_p, && \xi(x) \in \mathcal{D}_{\xi(x)}, t \in \mathcal{D}_{t > 0} \\
        &\frac{\partial y_p(t, x,\xi(x))}{\partial x}\Big|_{x = 500} = 0, && \xi(x) \in \mathcal{D}_{\xi(x)}, t \in \mathcal{D}_{t} \\
        & y_\theta(0, x, \xi(x)) = 1, && \xi(x) \in \mathcal{D}_{\xi(x)}, x \in \mathcal{D}_{x}
    \end{aligned}
    \label{eq:ald_boundary_conds}
\end{equation}
which enforce that the substrate not have any precursor initially, the ambient gaseous density be constant at $z_p \in \mathbb{R}$, no diffusion occurs beyond a $500 \mu m$ depth, and all the sites are initially available, respectively. We simulate \eqref{eq:ald_pdes} and \eqref{eq:ald_boundary_conds} with $z_p = 0.07$ in \texttt{InfiniteOpt.jl} that is evaluated via direct transcription over 10 random field samples, 100 spatial support points, and 10 temporal support points. The resulting finite model is solved using \texttt{KNITRO v12.4}. Figure \ref{fig:ald_sim} shows the final coverage profile $1-y_\theta(t, x, \xi(x))$ whose behavior is consistent with the results presented in \cite{keuter2015modeling}.

To setup the RFO problem, we seek to optimally choose $z_p$ such that the coverage profile $1-y_\theta(10, x, \xi(x))$ follows a desired setpoint $1-\bar{y}_\theta(x)$. Here the setpoint function is defined via a modified inverse sigmoid function that resembles the desired profile presented in \cite{keuter2015modeling}. 
\begin{equation}
    \label{eq: ald_setpoint}
    \bar{y}_\theta(x) = 1 - \frac{1}{1+\exp{(x - \frac{\rho_{1}}{2})}^{\rho_{2}}}
\end{equation}
where $\rho_{1}$ = 500 and $\rho_{2}$ = 0.15. The objective seeks to minimize the expected setpoint tracking error at the final time:
\begin{equation}
    \min \mathbb{E}_{xi(x)}\left[\int_{x \in \mathcal{D}_x} (y_\theta(10, x, \xi(x)) - \bar{y}_\theta(x))^2 dx \right].
    \label{eq:ald_objective}
\end{equation}
Putting together \eqref{eq:ald_pdes}, \eqref{eq:ald_boundary_conds}, and \eqref{eq:ald_objective} we obtain the RFO problem:
\begin{equation}
	\begin{aligned}
		&&\min_{} &&& \mathbb{E}_{xi(x)}\left[\int_{x \in \mathcal{D}_x} (y_\theta(10, x, \xi(x)) - \bar{y}_\theta(x))^2 dx \right] \\
		&& \text{s.t.} &&& \frac{\partial y_p(t, x, \xi(x))}{\partial t} = D \frac{\partial^{2} y_p(t, x, \xi(x))}{\partial x^2} - \gamma \xi(x) y_p(t, x, \xi(x)) y_\theta(t, x, \xi(x)), && \xi(x) \in \mathcal{D}_{\xi(x)}, (t,x) \in \mathcal{D}_{t,x} \\
        &&&&&\frac{\partial y_\theta(t, x, \xi(x))}{\partial t} = -\eta \xi(x) y_p(t, x, \xi(x)) y_\theta(t, x, \xi(x)), && \xi(x) \in \mathcal{D}_{\xi(x)}, (t,x) \in \mathcal{D}_{t,x} \\
        &&&&&y_p(0, x,\xi(x)) = 0, && \xi(x) \in \mathcal{D}_{\xi(x)}, x \in \mathcal{D}_{x} \\
        &&&&&y_p(t, 0,\xi(0)) = z_p, && \xi(x) \in \mathcal{D}_{\xi(x)}, t \in \mathcal{D}_{t > 0} \\
        &&&&&\frac{\partial y_p(t, x,\xi(x))}{\partial x}\Big|_{x = 500} = 0, && \xi(x) \in \mathcal{D}_{\xi(x)}, t \in \mathcal{D}_{t} \\
        &&&&& y_\theta(0, x, \xi(x)) = 1, && \xi(x) \in \mathcal{D}_{\xi(x)}, x \in \mathcal{D}_{x}.
	\end{aligned}
	\label{eq:ald_form}
\end{equation}
This is implemented in \texttt{InfiniteOpt.jl} using the same solution configuration used for simulation. We solve the deterministic variant of \eqref{eq:ald_form} that uses $\mu(x)$ for $\xi(x)$ and we simulate the system response with different realizations of $\xi(x)$ with the deterministic optimal value of $z_p$

\begin{figure}[!htb]
	\includegraphics[width=0.55\textwidth]{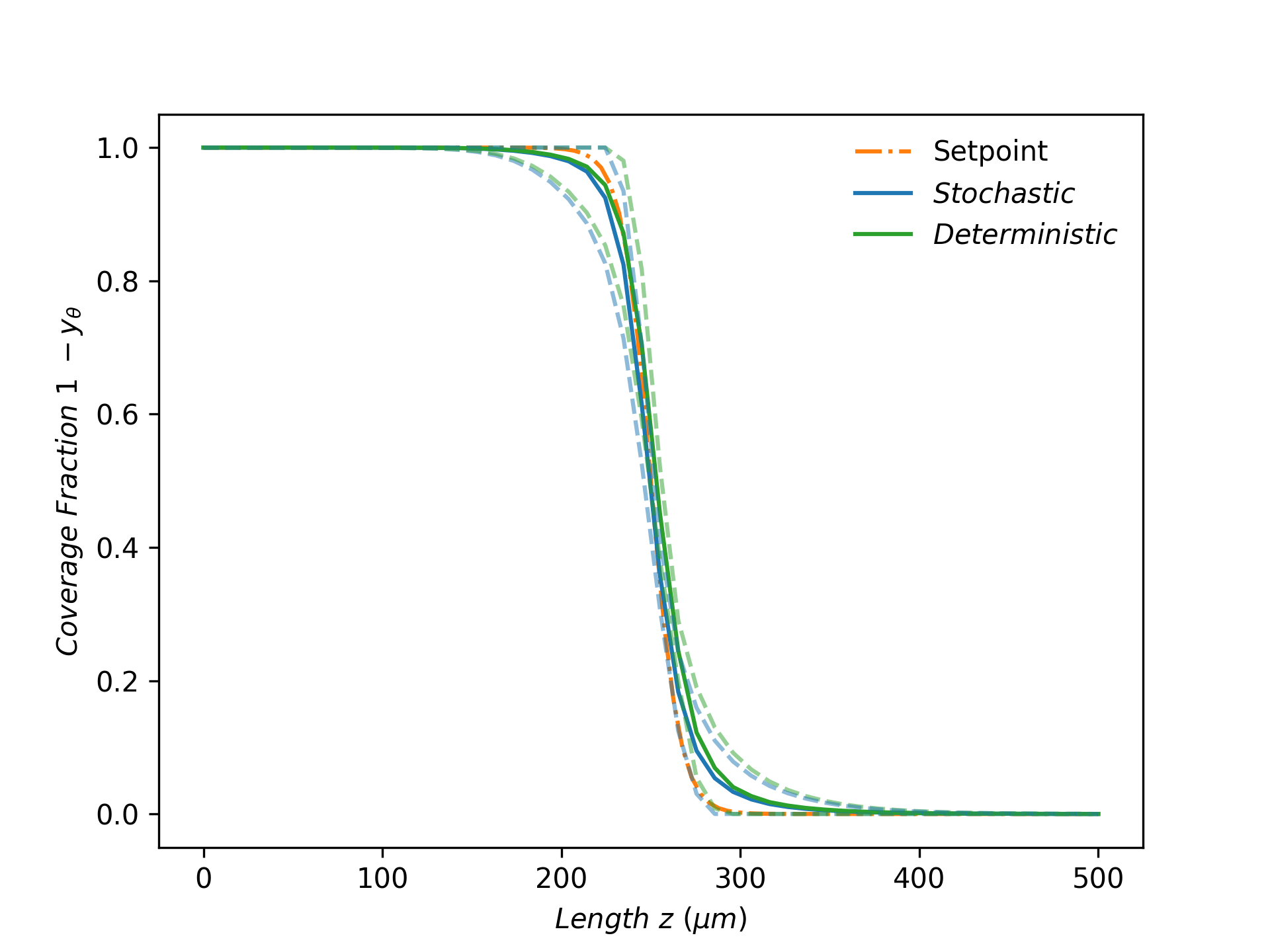}
	\centering
	\caption{Simulated ensemble coverage profiles utilizing the stochastic and deterministic solutions of $z_p$. Both are able to well track the desired setpoint.}
	\label{fig:ald_comparison}
\end{figure}

Figure \ref{fig:ald_comparison} shows the optimal coverage MC ensemble profiles in juxtaposition to the set-point profile. In both cases, the system is able to closely track the desired set-point profile. Moreover, we observe the deterministic and stochastic solutions to be very similar in this case since their choices of $z_p$ are nearly identical (the values are 0.0391 and 0.0393, respectively). Thus, for the behavior of the random field uncertainty $\xi(x)$ in this case, the RFO formulation is not more advantageous than the deterministic one. This highlights that accounting for the uncertainty via an RFO formulation will not always significantly outperform its deterministic variant. However, we can envision other random fields that would incur a significant difference. Thus, in this case study we have demonstrated how our RFO framework readily enables us to characterize uncertainty that propagates over spatial position, simulate a system that is subjected to this uncertainty, and then form an RFO problem to optimally choose conditions at which to run the process.

\subsection{Excursion Probability Optimization} \label{sec:diffusion_case}

\begin{figure}[!htb]
	\includegraphics[width=0.4\textwidth]{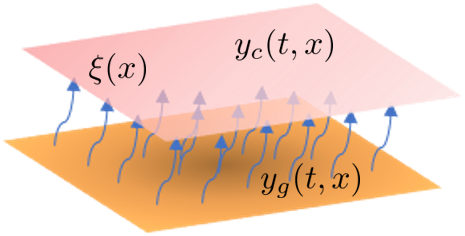}
	\centering
	\caption{Representative schematic of the transient diffusion system considered in Section \ref{sec:diffusion_case}.}
	\label{fig:diffusion_schematic}
\end{figure}

We further exemplify the principles discussed above by implementing an RFO problem that seeks to minimize the excursion probability of violating a certain concentration/temperature threshold (we refer to it in terms of temperature) in a diffusion system. In particular, we consider a 2D transient diffusion system with temperature $y_c(t, x, \xi(x)) \in \mathbb{R}_+$ subject to spatially random diffusivity $\xi(x) \in \mathcal{D}_{\xi(x)}$ (a 2D random field) that we seek to heat with a spatially variable heating plate with temperature $y_g(t, x) \in [0, 0.1]$. Figure \ref{fig:diffusion_schematic} provides an illustrative schematic of this system.

\begin{figure}[!htb]
	\includegraphics[width=0.8\textwidth]{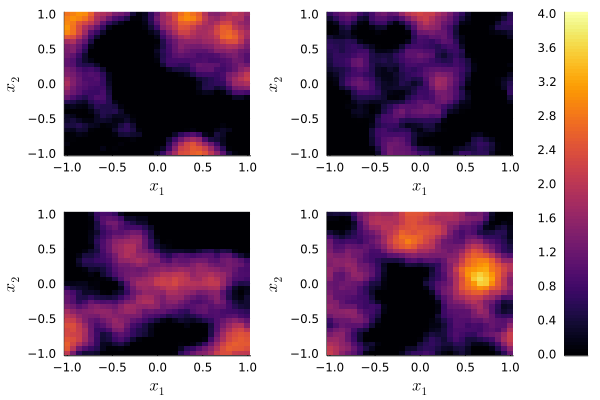}
	\centering
	\caption{Realizations of the random diffusivity field $\hat{\xi}_k(x)$.}
	\label{fig:diffusion_fields}
\end{figure}

The random field diffusivity is characterized as a Gaussian random field with moments:
\begin{equation}
	\begin{aligned}
		&\mu(x) = 0.5 \\
		&\Sigma_\text{M}(x, x') = \frac{2^{1-\nu}}{\Gamma(\nu)}\left(\frac{\sqrt{2\nu}||x-x'||}{\beta}\right)^\nu K_\nu\left(\frac{\sqrt{2\nu}||x-x'||}{\beta}\right)
	\end{aligned}
	\label{eq:diffusion_moments}
\end{equation}
where we set $\beta = 0.25$ and $\nu = 1.5$. We depict four sample function realizations of $\xi(x)$ in Figure \ref{fig:diffusion_fields}.

Drawing from the transient diffusion relations in \eqref{eq:transient_diffusion}, we model our system: 
\begin{equation}
	\frac{\partial y_c(t, x, \xi(x))}{\partial t} = \xi(x) \nabla^2_x y_c(t, x, \xi(x)) + y_g(t, x), \ \xi(x) \in \mathcal{D}_{\xi(x)}, (t, x) \in \mathcal{D}_{t,x}
	\label{eq:case_diffusion}
\end{equation}
where $\mathcal{D}_{t,x} = \mathcal{D}_t \times \mathcal{D}_x = [0, 1] \times [-1, 1]^2$. For boundary conditions, we use Dirichlet conditions:
\begin{equation}
	\begin{aligned}
		&y_c(0, x, \xi(x)) = 0, && \xi(x) \in \mathcal{D}_{\xi(x)}, x \in \mathcal{D}_x\\
		&y_c(t, [-1, x_2], \xi([-1, x_2])) = 0, && \xi(x) \in \mathcal{D}_{\xi(x)}, t \in \mathcal{D}_t, x_2 \in \mathcal{D}_{x_2} \\
		&y_c(t, [1, x_2], \xi([1, x_2])) = 0, && \xi(x) \in \mathcal{D}_{\xi(x)}, t \in \mathcal{D}_t, x_2 \in \mathcal{D}_{x_2} \\
		&y_c(t, [x_1, -1], \xi([x_1, -1])) = 0, && \xi(x) \in \mathcal{D}_{\xi(x)}, t \in \mathcal{D}_t, x_1 \in \mathcal{D}_{x_1} \\
		&y_c(t, [x_1, 1], \xi([x_1, 1])) = 0, && \xi(x) \in \mathcal{D}_{\xi(x)}, t \in \mathcal{D}_t, x_1 \in \mathcal{D}_{x_1}.
	\end{aligned}
	\label{eq:diffusion_conditions}
\end{equation}
We seek to optimally raise and maintain the plate temperature to a setpoint $y_c(t, x, \xi(x)) = 0.2$ by enforcing the objective:
\begin{equation}
	\min \mathbb{E}_{\xi(x)}\left[\int_{(t,x) \in \mathcal{D}_{t,x}}(y_c(t,x) - 0.2)^2 dt\right].
	\label{eq:diffusion_objective}
\end{equation}
Moreover, we limit how the plate temperature exceeds a safety threshold $u = 0.25$ via the excursion probability constraint:
\begin{equation}
	\mathbb{P}_{\xi(x)}\left(\max_{(t, x) \in \mathcal{D}_{t,x}} y_c(t,x, \xi(x)) > 0.25 \right) \leq \alpha 
	\label{eq:diffusion_excursion}
\end{equation}
where $\alpha \in [0, 1]$ is a probability level. Putting together \eqref{eq:case_diffusion} - \eqref{eq:diffusion_excursion}, we obtain the RFO problem:
\begin{equation}
	\begin{aligned}
		&&\min_{y_c(\cdot), y_g(\cdot)} &&&\mathbb{E}_{\xi(x)}\left[\int_{(t,x) \in \mathcal{D}_{t,x}}(y_c(t,x) - 0.2)^2 dt\right] \\
		&& \text{s.t.} &&& \frac{\partial y_c(t, x, \xi(x))}{\partial t} = \xi(x) \nabla^2_x y_c(t, x, \xi(x)) + y_g(t, x), && \xi(x) \in \mathcal{D}_{\xi(x)}, (t, x) \in \mathcal{D}_{t,x} \\
		&&&&&y_c(0, x, \xi(x)) = 0, && \xi(x) \in \mathcal{D}_{\xi(x)}, x \in \mathcal{D}_x\\
		&&&&&y_c(t, [-1, x_2], \xi([-1, x_2])) = 0, && \xi(x) \in \mathcal{D}_{\xi(x)}, t \in \mathcal{D}_t, x_2 \in \mathcal{D}_{x_2} \\
		&&&&&y_c(t, [1, x_2], \xi([1, x_2])) = 0, && \xi(x) \in \mathcal{D}_{\xi(x)}, t \in \mathcal{D}_t, x_2 \in \mathcal{D}_{x_2} \\
		&&&&&y_c(t, [x_1, -1], \xi([x_1, -1])) = 0, && \xi(x) \in \mathcal{D}_{\xi(x)}, t \in \mathcal{D}_t, x_1 \in \mathcal{D}_{x_1} \\
		&&&&&y_c(t, [x_1, 1], \xi([x_1, 1])) = 0, && \xi(x) \in \mathcal{D}_{\xi(x)}, t \in \mathcal{D}_t, x_1 \in \mathcal{D}_{x_1} \\
		&&&&& \mathbb{P}_{\xi(x)}\left(\max_{(t, x) \in \mathcal{D}_{t,x}} y_c(t,x, \xi(x)) > 0.25 \right) \leq \alpha 
	\end{aligned}
	\label{eq:diffusion_rfo}
\end{equation}
which seeks a tradeoff policy $y_g^*(t,x)$ that optimally minimizes the expected tracking error while not exceeding the excursion threshold. This inherent bi-objective trade-off behavior is controlled with the choice of $\alpha$. Following the approach in \cite{pulsipher2019scalable}, we invert Problem \eqref{eq:diffusion_rfo} to obtain the equivalent $\epsilon$-constrained form:
\begin{equation}
	\begin{aligned}
		&&\min_{y_c(\cdot), y_g(\cdot)} &&& \mathbb{P}_{\xi(x)}\left(\max_{(t, x) \in \mathcal{D}_{t,x}} y_c(t,x, \xi(x)) > 0.25 \right) \\
		&& \text{s.t.} &&& \frac{\partial y_c(t, x, \xi(x))}{\partial t} = \xi(x) \nabla^2_x y_c(t, x, \xi(x)) + y_g(t, x), && \xi(x) \in \mathcal{D}_{\xi(x)}, (t, x) \in \mathcal{D}_{t,x} \\
		&&&&&y_c(0, x, \xi(x)) = 0, && \xi(x) \in \mathcal{D}_{\xi(x)}, x \in \mathcal{D}_x\\
		&&&&&y_c(t, [-1, x_2], \xi([-1, x_2])) = 0, && \xi(x) \in \mathcal{D}_{\xi(x)}, t \in \mathcal{D}_t, x_2 \in \mathcal{D}_{x_2} \\
		&&&&&y_c(t, [1, x_2], \xi([1, x_2])) = 0, && \xi(x) \in \mathcal{D}_{\xi(x)}, t \in \mathcal{D}_t, x_2 \in \mathcal{D}_{x_2} \\
		&&&&&y_c(t, [x_1, -1], \xi([x_1, -1])) = 0, && \xi(x) \in \mathcal{D}_{\xi(x)}, t \in \mathcal{D}_t, x_1 \in \mathcal{D}_{x_1} \\
		&&&&&y_c(t, [x_1, 1], \xi([x_1, 1])) = 0, && \xi(x) \in \mathcal{D}_{\xi(x)}, t \in \mathcal{D}_t, x_1 \in \mathcal{D}_{x_1} \\
		&&&&& \mathbb{E}_{\xi(x)}\left[\int_{(t,x) \in \mathcal{D}_{t,x}}(y_c(t,x) - 0.2)^2 dt\right] \leq \epsilon
	\end{aligned}
	\label{eq:diffusion_inverted}
\end{equation}
where $\epsilon \in \mathbb{R}_+$ is the tracking error budget (i.e., Pareto parameter). We can reformulate the excursion probability following the big-$M$ constraint approach discussed in Section \ref{sec:direct_transcription}:
\begin{equation}
	\mathbb{P}_{\xi(x)}\left(\max_{(t, x) \in \mathcal{D}_{t,x}} y_c(t,x, \xi(x)) > 0.25 \right) = \mathbb{E}_{\xi(x)}\left[y_b(\xi(x))\right]
\end{equation}
with constraints:
\begin{equation}
	\begin{aligned}
		&\overline{y}_c(\xi(x)) \geq y_c(t, x, \xi(x)), && \xi(x) \in \mathcal{D}_{\xi(x)}, (t, x) \in \mathcal{D}_{t,x} \\
		&\overline{y}_c(\xi(x)) - 0.25 \leq y_b(\xi(x))M, && \xi(x) \in \mathcal{D}_{\xi(x)}
	\end{aligned}
\end{equation}
where $y_b(\xi(x)) \in \{0, 1\}$, $\overline{y}_c(\xi(x)) \in \mathbb{R}_+$, and $M \in \mathbb{R}_+$ is a suitable upper bound. This reformulation introduces the binary variable $y_b(\xi(x))$ which compounds the problem complexity, but this can be alleviated via the continuous relaxation proposed in \cite{pulsipher2019scalable} and discussed in Section \ref{sec:direct_transcription} where we use the relaxation $y_b(\xi(x)) \in [0, 1]$ to obtain Pareto solutions which we can then round to integrality using an appropriate rounding rule (often recovering the optimal mixed-integer solution to high accuracy).
\\

We implement \eqref{eq:diffusion_inverted} using the continuous big-$M$ constraint reformulation of the excursion probability in \texttt{InfiniteOpt.jl} to obtain the Pareto pairs that correspond to: 
\begin{equation}
	\epsilon \in \{0.0947,0.095,0.096,0.097,0.098\}.
	\label{eq:epsilon_values}
\end{equation}
We transform the model via direct transcription using 7 realizations of $\xi(x)$ and 9,600 grid points over $(t,x) \in \mathcal{D}_{t,x}$ (10 for $t$ and 31 for $x_1,x_2$), and is solved using \texttt{Gurobi v.9.1.2}. The coarseness of the transcription grid shows how the combinatorics behind such transformations quickly increase the problem size. Hence, more sophisticated transformation methods are required for higher fidelity transformations as discussed in Section \ref{sec:rfo_other_transforms}. 

\begin{figure}[!htb]
	\includegraphics[width=0.6\textwidth]{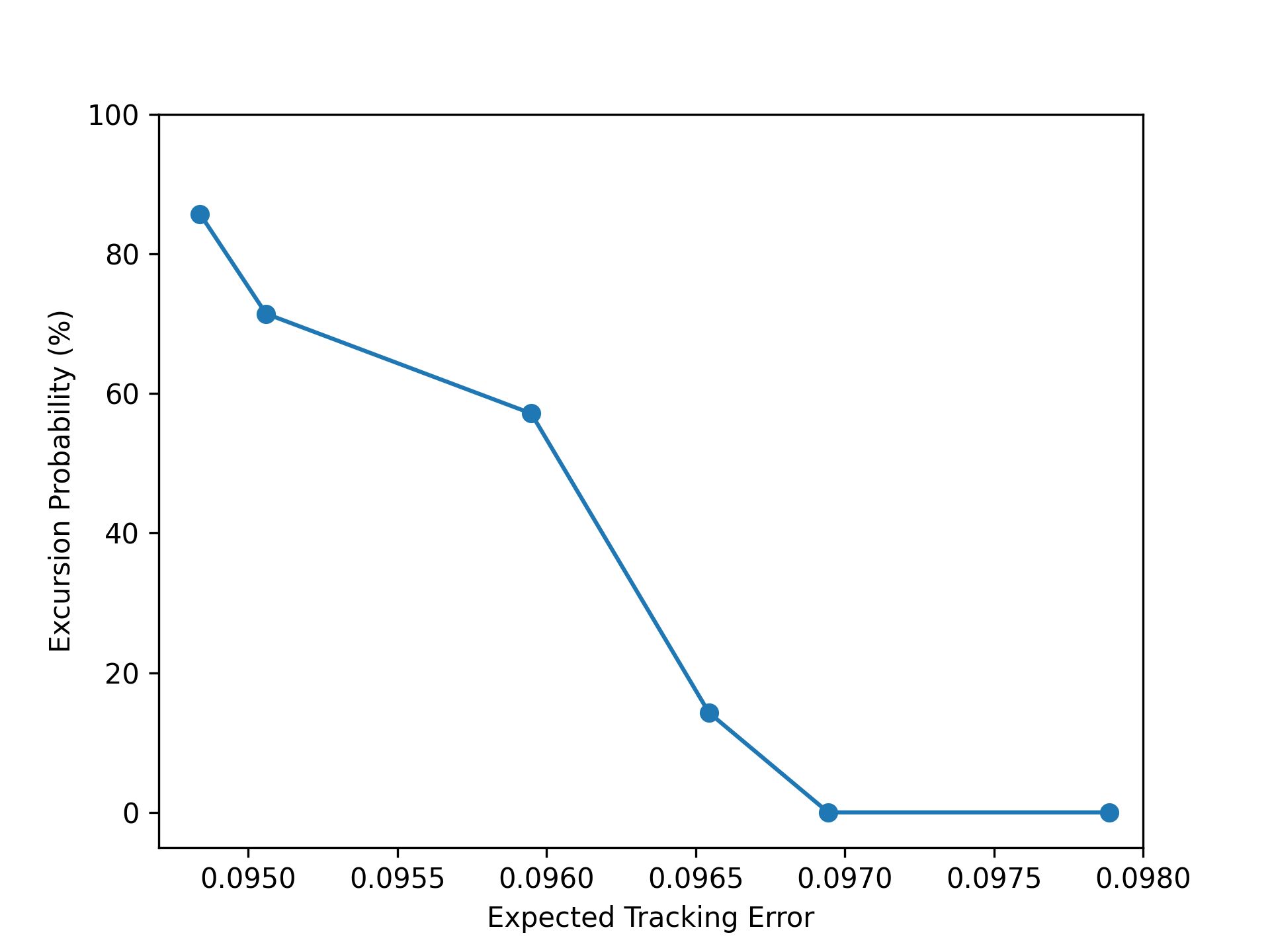}
	\centering
	\caption{Pareto frontier capturing tradeoff between minimizing the excursion probaility in \eqref{eq:diffusion_excursion} and the expected tracking error of \eqref{eq:diffusion_objective}.}
	\label{fig:diffusion_pareto}
\end{figure}

\begin{figure}[!htb]
    \centering
	\begin{subfigure}[b]{0.45\textwidth}
        \centering
        \includegraphics[width=0.8\textwidth]{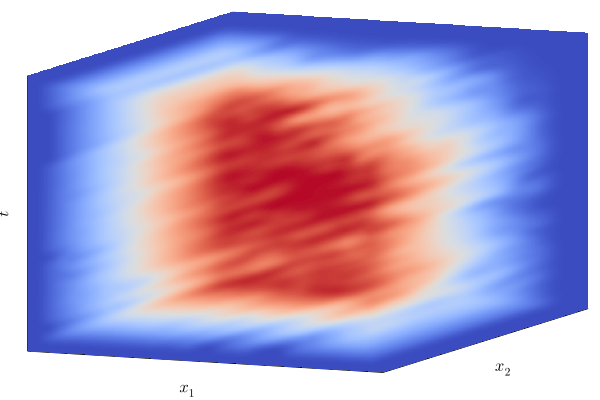}
        \caption{$y_g^*(t, x)$}
    \end{subfigure}
	\quad
    \begin{subfigure}[b]{0.45\textwidth}
        \centering
        \includegraphics[width=0.8\textwidth]{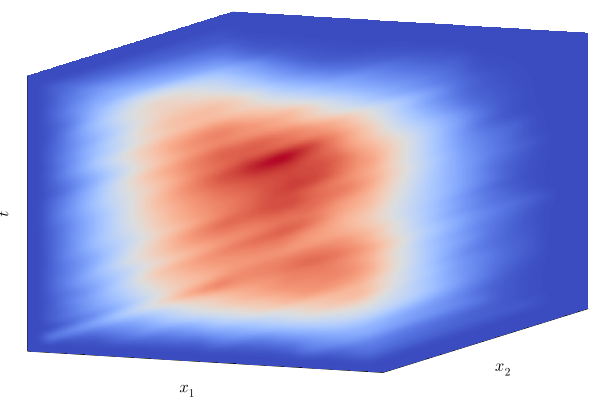}
        \caption{$\mathbb{E}_{\xi(x)}[y_c^*(t, x, \xi(x))]$}
    \end{subfigure}
    \caption{Optimal heating policy and response for problem \eqref{eq:diffusion_inverted},  corresponding to $\epsilon = 0.098$.}
    \label{fig:diffusion_heatmaps}
\end{figure}

Figure \ref{fig:diffusion_pareto} shows the Pareto frontier obtained from the values of $\epsilon$ in \eqref{eq:epsilon_values}. This clearly exhibits the tradeoff between tightly tracking the temperature setpoint and minimizing the temperature excursion relative to the sample function realizations of $\xi(x)$ used in this study. Figure \ref{fig:diffusion_heatmaps} shows the optimal Pareto pair corresponding to $\epsilon = 0.098$ where we observe the compromised deterministic heating policy $y_g^*(t, x)$ that is used to obtain the optimal response $y_c^*(t,x, \xi(x))$. This highlights how we can use excursion probabilities to shape the policies derived from RFO problems in a manner analogous to what is accomplished with chance constraints in stochastic programming.

\section{Conclusions and Future Work} \label{sec:conclusion}

We have presented a general modeling framework for leveraging random field uncertainty in an optimization context that allows us to conduct uncertainty propagation over general infinite domains (e.g., space and/or time). This framework, which we call random field optimization, provides a unified approach to model and optimize engineering systems that are affected by complex sources of uncertainty (e.g., stochastic processes or spatial random fields). Moreover, our proposed framework for modeling and optimizing RFO problems is general, capturing existing approaches such as multi-stage stochastic optimization as special cases. Our hope is that this framework will serve as a foundation for further advancement and innovation for tackling this challenging problem class.
\\

In future work, we plan to develop more advanced transformation approaches to achieve higher fidelity solutions in a more tractable manner. Generalizing the decomposition techniques proposed in \cite{shin2019scalable} and the order-reduction techniques common to PDE-constrained optimization provide promising avenues of research in this respect. We also plan to more fully implement this framework and these proposed transformation techniques in \texttt{InfiniteOpt.jl} in an effort to make RFO problems more accessible for future research and non-expert practitioners.

\section*{Acknowledgments}
This work was supported by the U.S. Department of Energy under grant DE-SC0014114.


\bibliography{references}
\end{document}